\documentclass[twoside, 11pt]{article}

\usepackage{amssymb, amsmath, latexsym, mathrsfs, verbatim, calc}
\usepackage{color}

\numberwithin{equation}{section}
\newtheorem{thm}{Theorem}[section]
\newtheorem{lem}[thm]{Lemma}

\newtheorem{prop}[thm]{Proposition}
\newtheorem{rem}[thm]{Remark}

\newtheorem{defn}[thm]{Definition}
\newtheorem{assum}[thm]{Assumption}

\definecolor{blue-violet}{rgb}{0.54, 0.17, 0.89}
\definecolor{purple1}{rgb}{0.63, 0.36, 0.94}
\definecolor{purple2}{rgb}{0.87, 0.0, 1.0}
\definecolor{purple}{rgb}{0.5, 0.0, 0.5}
\definecolor{pansypurple}{rgb}{0.47, 0.09, 0.29}
\definecolor{orange}{rgb}{1.0, 0.27, 0.0}

\def\ba{\begin{array}}
\def\ea{\end{array}}
\def\beq{\begin{equation}}
\def\bes{\begin{equation*}}
\def\ees{\end{equation*}}
\def\bea{\begin{eqnarray}}
\def\eea{\end{eqnarray}}
\def\beaa{\begin{eqnarray*}}
\def\eeaa{\end{eqnarray*}}

\def\dis{\displaystyle}

\def\no{\noindent}

\def\lastline{\par \vspace{-7.3ex} \no}

\def\nts{\negthinspace}

\def\ms{\medskip}
\def\bs{\bigskip}
\def\q{\quad}
\def\qq{\qquad}

\def\ol{\overline}
\def\ul{\underline}

\def\da{\mathop{\downarrow}}

\def\={=\nts \nts=\nts \nts=\nts \nts=}

\def\({\textnormal{(}}
\def\){\textnormal{)}}

\def\cd{\cdot}
\def\cds{\cdots}

\def\pa{\partial}

\def\a{\alpha}
\def\g{\gamma}
\def\d{\delta}
\def\e{\varepsilon}

\def\l{\lambda}
\def\m{\mu}
\def\n{\nu}

\def\si{\sigma}
\def\t{\tau}
\def\f{\varphi}
\def\th{\theta}
\def\o{\omega}

\def\f{\phi}
\def\vf{\varphi}
\def\vsi{\varsigma}
\def\p{\psi}

\def\D{\Delta}

\def\O{\Omega}

\def\cF{{\cal F}}
\def\cG{{\cal G}}

\def\cO{{\cal O}}

\def\cS{{\cal S}}

\def\hC{\mathbb{C}}

\def\hE{\mathbb{E}}
\def\hF{\mathbb{F}}

\def\hN{\mathbb{N}}

\def\hP{\mathbb{P}}

\def\hR{\mathbb{R}}

\def\hX{\mathbb{X}}

\def\sA{\mathscr{A}}
\def\sB{\mathscr{B}}
\def\sC{\mathscr{C}}
\def\sD{\mathscr{D}}

\def\sH{\mathscr{H}}

\def\sL{\mathscr{L}}

\def\sP{\mathscr{P}}

\def\sU{\mathscr{U}}

\def\liminf{\mathop{\ul{\rm lim}}}
\def\limsup{\mathop{\ol{\rm lim}}}

\def\neg{\negthinspace}

\def\no{\noindent}

\def\ms{\medskip}
\def\bs{\bigskip}
\def\q{\quad}
\def\qq{\qquad}

\def\pa{\partial}
\def\cd{\cdot}
\def\cds{\cdots}

\def\qed{\hfill \rule[0cm]{.25cm}{.25cm}\medskip}   
\def\dfnn{\stackrel{\triangle}{=}}

\def\b1{{\bf 1}}

\newenvironment{itm}{\vspace{-1ex}\begin{itemize}}{\end{itemize}}
\def\bi{\begin{itm}}
\def\ei{\end{itm}}

\def\equ_ind{\arabic{section}.\arabic{equation}}
\def\sec_ind{\arabic{section}}

\begin{document}
\title{\bf Optimal Dividend and Investment Problems under Sparre Andersen Model}

\author{Lihua Bai\thanks{\noindent School of Mathematical Sciences,
Nankai University, Tianjin, China. E-mail:
lhbai@nankai.edu.cn. This author is supported in part by Chinese NSF
grants \#11471171 and \#11571189.}, ~Jin Ma\thanks{\noindent
 Department of Mathematics,
University of Southern California, Los Angeles, CA 90089. E-mail:
jinma@usc.edu. This author is supported in part by US NSF grant  \#1106853. }, ~  Xiaojing Xing
\thanks{\noindent Department of Mathematics, University of Southern California, Los
Angeles, CA 90089. E-mail: xiaojinx@usc.edu.
}}

\date{}
\maketitle

\vspace{-2mm}
\begin{abstract}
In this paper we study a class of optimal dividend and investment
problems assuming that the underlying reserve process follows the Sparre Andersen model, that is,
the claim frequency is a ``renewal" process, rather than a standard compound Poisson process. The main feature of such  problems is that the underlying reserve dynamics, even in its simplest form, is no longer Markovian.
By using the  {\it backward Markovization} technique we recast the problem in a Markovian framework with
expanded dimension representing the time elapsed after the last claim, with which we investigate the regularity of the value function, and validate the dynamic programming principle. Furthermore, we show that the value function
is the unique {\it constrained viscosity solution} to the associated HJB equation on a cylindrical domain
on which the problem is well-defined.
\end{abstract}

\vfill \bs

\no

{\bf Keywords.} \rm Optimal dividend problem, Sparrer Andersen model, backward Markovization,
dynamic programming,  Hamilton-Jacobi-Bellman equation, constrained viscosity solution.

\bs

\no{\it 2000 AMS Mathematics subject classification:} 91B30, 93E20;
60K05, 35D40.

\eject
\section{Introduction}
\setcounter{equation}{0}

The problem of maximizing the cumulative discounted dividend payout can be traced back to the seminal work of de Finetti \cite{Fen} in 1957, when he proposed to measure the performance of an insurance portfolio by looking at
the maximum possible dividend paid during its lifetime, instead of focussing only on the safety aspect measured by its ruin probability.
Although other criteria such as the  so-called Gordon model \cite{GOR} as well as the simpler model by Miller-Modigliani \cite{MIL} have been proposed over the years, to date  the cumulative discounted dividend is still widely accepted as an important
and useful performance index, and various approaches have been employed to find the optimal strategy that maximizes
such index.  The solution of the optimal dividend problem under the classical Cram\'er-Lundberg model has been obtained in various forms. Gerber \cite{Ger} first
showed that an optimal dividend strategy has a ``band" structure. Since then the optimal dividend policies, especially the barrier strategies, have been investigated in various settings, sometimes under more  general reserve models (see, e.g.,  \cite{albrh,alb, AHT, Ger2, HoTak2, jea, loe,mni,sch1},  to mention
 a few).
 We refer the interested reader to the excellent 2009 survey by  Albrecher-Thonhauser \cite{AT} and the exhaustive references cited therein for the past developments on this issue.

The more general optimization problems for insurance models involving the possibility of investment and/or reinsurance
have also been studied quite extensively in the past two decades.
In 1995, Browne \cite{bro} first considered the problem of minimizing  the probability of ruin under a diffusion approximated Cram\'er-Lundberg model, where the insurer is allowed to invest some fraction of the reserve dynamically into a Black-Scholes market.
Hipp-Plum \cite{hip} later considered the same problem with a compound Poisson claim process. The problems involving either proportional or excess-of-loss reinsurance strategies have also been been studied under the Cram\'er-Lundberg model or its
diffusion approximations (see, e.g., \cite{hip1, HoTak1, HoTak2, sch}).
%
The optimal dividend and reinsurance problem with transaction cost and taxes was studied by the first author of this paper with various co-authorships \cite{BGZ, BMP, BP}; whereas the ruin problems, reinsurance problems, and universal variable
insurance problems involving investment in the more general jump diffusion framework have been investigated by the second author \cite{LiuMa, MS, MYu}, from the stochastic control perspective.
We should remark that the two references that are closest to the present paper are Azcue-Muler \cite{azc, azcu}, obtained in 2005 and
2010, respectively. The former concerns the optimal dividend-reinsurance, and the latter
concerns the optimal dividend-investment. Both papers followed the dynamic programming approach, and the analytic properties
of the value function, including its being the viscosity solution to  the associated Hamilton-Jacobi-Bellman (HJB) equation became
the main purpose.


It is worth noting, however, that all  aforementioned results are based on  the Cram\'er-Lundberg type of
surplus dynamics or its  variations within the Markovian paradigm, whose analytical structure plays a fundamental role.
A well-recognized  generalization of such model is one in which the Poisson claim number process is replaced by a {\it renewal} process, known as the Sparre Andersen risk model \cite{spa}. The dividend problem under such a
model is much subtler due to its {\it non-Markovian} nature in general, and the literature is much more limited. In this context,
Li-Garrido \cite{li} first studied the properties of the
renewal risk reserve process with a barrier strategy. Later, after  calculating the moments of the expected discounted dividend payments under a barrier strategy in \cite{albr},  Albrecher-Hartinger \cite{albrh} showed that, unlike the classical Cram\'er-Lundberg model,  even in the case of Erlang(2) distributed interclaim times and exponentially distributed claim amounts,
the horizontal barrier strategy is no longer optimal. Consequently, the optimal dividend problem under the Sparre Andersen
models has since been listed as an open problem that requires attention (see \cite{AT}), and to the best of our knowledge, it remains unsolved to this day.

The main technical difficulties, from the stochastic control perspective, for a general optimal dividend problem under Sparre Andersen model can be roughly summarized into two major points: the non-Markovian nature of the model, and the random duration of the insurance portfolio. We note that although  the former would seemingly invalidate the dynamical programming approach, a Markovization is possible, by extending the dimension of the state space of the risk process, taking into account the time elapsed since the last claim (see \cite{AT}). It turns out that such an extra variable would cause some subtle technical difficulties in analyzing the regularity of the value function. For example, as we shall see later, unlike the compound Poisson cases studied in \cite{azc, azcu}, even the continuity of the value function requires some heavy arguments, much less the Lipschitz properties which plays
a fundamental role in a standard argument. For the latter issue, since we are
focusing on the life of the portfolio until ruin, the optimization problem naturally has a random terminal time. While it is known in theory that such a problem can often be converted to one with a fixed (deteministic)
terminal time (see, e.g., \cite{BEML2008}) once the distribution of the random terminal is known, finding the distribution for
the ruin time under Sparre Andersen model is itself a challenging problem, even under very explicit strategies (see, e.g., \cite{albr, Ger2, li}), which makes the optimization problem technically prohibitive along this line.

This paper is our first attempt to attack this open problem. We will start with a  rather simplified renewal reserve model but
allowing both investment and dividend payments. As was suggested in \cite{AT}, our plan is to first ``Markovize" the model and then study the optimal dividend problem via the dynamic programming approach. Specifically, we shall first investigate the property of the value function and then validate the  {\it dynamic programming principle} (DPP), from which we can formally derive the associated HJB equation to which the value function is a solution  in some sense. An important observation, however, is that the value function could very well be discontinuous at
the boundary of a region on which it is well-defined, and no explicit boundary condition can be established directly from the information of the problem. Among other things, the lack of  boundary information of the HJB equation will make the comparison principle, whence uniqueness, particularly subtle, if not impossible.
To overcome this difficulty we shall invoke the notion of {\it constrained viscosity solution} for the exit problems (see, e.g.,
Soner \cite{Soner}), and as it turns out we can prove that the value function is indeed a constrained viscosity solution to the associated HJB equation
on an appropriately defined domain, completing the dynamic programming approach on this problem. To the best of
our knowledge, these results are novel.

%

The rest of the paper is organized as follows. In section 2, we establish the basic setting, formulate the problem, and introduce the backward Markovization technique. In section 3 we study the properties of the value function and prove the continuity of
the value function in the temporal variable. In Sections 4 and 5 we prove the continuity of the value function in spacial variable $x$ and the delayed variable $w$, respectively. In Section 6  we validate the Dynamic Programming Principle (DPP), and in
Section 7 we show that the value function is a constrained viscosity solution to the HJB equation. Finally, in section 8 we prove
the comparison principle, hence prove that the value function is the  unique constrained viscosity solution among a fairly
general class of functions.

%

\section{Preliminaries and Problem Formulation}
\setcounter{equation}{0}

Throughout this paper we assume that all uncertainties come from
a common complete probability space $(\O,\cF, \hP)$ on which is
defined  $d$-dimensional Brownian motion $B=\{B_t:t\ge 0\}$, and a {\it renewal}
counting process  $N=\{N_t\}_{t\ge0}$, independent of $B$. More precisely, we denote
$\{\sigma_{n}\}_{n=1}^\infty$ to be the jump times ($\sigma_{0}:=0$) of the counting
process $N$, and $T_{i}=\sigma_{i}-\sigma_{i-1}$, $i=1,2,\cdots$ to be its waiting times
(the time elapses between successive jumps). We assume that $T_i$'s are independent
and identically distributed,
with a common distribution $F:\hR_+\mapsto \hR_+$; and that there exists an {\it intensity function} $\l:[0,\infty)\mapsto [0,\infty)$  such that
$\bar F(t):=\hP\{T_1>t\}=\exp\{-\int_0^t\lambda(u)du\}$. In other words, $\lambda(t)=f(t)/\bar F(t)$, $t\ge 0$, where $f$ is the common density
function of $T_i$'s.

Further, throughout the paper we will denote, for a  generic Euclidean space
$\hX$, regardless of its dimension, $(\cd,\cd)$ and $|\cd|$ to be its inner product and norm, respectively.
Let $T>0$ be a given time horizon, we denote the space of continuous functions taking values in $\hX$ with the usual sup-norm by $C([0,T];\hX)$, and we shall make use
of the following notations:
\begin{itemize}
\item{} For any sub-$\si$-field $\cG\subseteq\cF$ and $1\le p<\infty$,
$L^p(\cG;\hX)$ denotes the space of all $\hX$-valued, $\cG$-measurable
random variables $\xi$ such that $\hE|\xi|^p<\infty$. As usual, $\xi\in L^\infty
(\cG;\hX)$ means that it is a bounded, $\cG$-measurable random variable.

\item{} For a given filtration $\hF=\{\cF_t:t\ge0\}$ in $\cF$, and $1\le p<\infty$, $L^p_\hF([0,T];\hX)$ denotes the space of all
$\hX$-valued, $\hF$-progressively measurable processes $\xi$ satisfying
$\hE\int_0^T|\xi_t|^pdt<\infty$. The meaning of $L^\infty _\hF([0,T];\hX)$
is defined similarly.
\end{itemize}

\subsection{Backward Markovization and Delayed Renewal Process}

An important ingredient of the Sparre Andersen model, is the following ``compound renewal process"
that will be used to represent the {\it claim process} in our reserve mode: $Q_t=\sum_{i=1}^{N_t}U_i$, $t\ge 0$, where $N$ is the
renewal process representing the frequency of the incoming claims, whereas $\{U_{i}\}_{i=1}^\infty$ is a sequence of
random variables representing the ``severity" (or claim size) of the incoming claims. We assume that $\{U_i\}$  are independent, identically
distributed with a common distribution $G:\hR_+\mapsto \hR_+$, and are independent of $(N, B)$.

The main feature of the Sparre Andersen model, which fundamentally differentiate this paper with all existing works is that
the process $Q$  is non-Markovian
in general (unless the counting process $N$ is a Poisson process), consequently we can not directly apply
the dynamic programming approach. We shall therefore first apply
the so-called {\it Backward Markovization} technique (cf. e.g., 
\cite{RSST}) to overcome this obstacle. More precisely,
we define a new process
\bea
\label{W}
W_t=t-\sigma_{N_t}, \qq t\ge 0,
\eea
 be the time elapsed since the last jump. Then clearly, $0\le W_t\le t\le T$, for $t\in [0,T]$, and it is known (see, e.g., \cite{RSST}) that
the process $(t, Q_t, W_t)$, $t\ge0$, is a piecewise deterministic Markov process (PDMP). We note that at each jump time $\si_i$, the jump size $|\D W_{\si_i}|=\si_i-\si_{i-1}=T_i$.

Throughout this paper we shall consider the following filtration $\{\cF_t\}_{t\ge 0}$, where $\cF_t:=\cF^B_t\vee\cF^Q_t\vee\cF^W_t$, $t\ge 0$. Here $\{\cF^\xi_t:t\ge 0\}$ denotes the natural filtration generated by process $\xi=B, Q, W$, respectively, with the usual $\hP$-augmentation such that it satisfies the {\it
usual hypotheses} (cf. e.g.,  \cite{PR}).

A very important element in the study of the dynamic optimal control problem with final horizon is to allow the starting point to
be any time $s\in [0,T]$. In fact, this is one of the main subtleties in the Sparre Andersen model, which we now describe.
Suppose that, instead of starting the clock at $t=0$, we start from $s\in[0,T]$, such that $W_s=w$, $\hP$-a.s. Let us consider the
{\it regular conditional probability distribution} (RCPD) $\hP_{sw}(\cd):=\hP[\,\cd\,|W_s=w]$ on $(\O, \cF)$, and consider the
``shifted" version of processes $(B, Q, W)$ on the space $(\O, \cF, \hP_{sw}; \hF^s)$, where $\hF^s=\{\cF_t\}_{t\ge s}$.
We first define $B^s_t:=B_{t}-B_s$, $t\ge s$. Clearly, since $B$ is independent of $(Q, W)$,
$B^s$ is  an $\hF^s$-Brownian motion under $\hP_{sw}$, defined on $[s,T]$, with $B^s_s=0$. Next, we restart the clock at time $s\in[0,T]$ by defining the new counting process $N^s_t:=N_{t}-N_s$, $t\in [s,T]$.  Then, under $\hP_{sw}$, $N^s$ is a ``delayed" renewal process, in the sense that
while its waiting times $T^s_i$, $i\ge 2$, remain independent, identically distributed as the original $T_i$'s, its ``time-to-first jump", denoted by
$T^{s,w}_1:=T_{N_s+1}-w=\si_{N_s+1}-s$, should have the survival probability
\bea
\label{PTsw}
\hP_{sw}\{T^{s,w}_1>t\}=\hP\{T_1>t+w|T_1>w\}=e^{\int_w^{w+t}\l(u)du}.
\eea
 In what follows we shall denote $N^s_t\big|_{W_s=w}:=N^{s,w}_t$,
$t\ge s$, to emphasize the dependence on $w$ as well. Correspondingly, we shall denote $Q^{s,w}_t=\sum_{i=1}^{N^{s,w}_t}U_i$
and $W^{s,w}_t:=w+W_{t}-W_s=w+[(t-s)-(\si_{N_{t}}-\si_{N_s})]$, $t\ge s$. It is readily seen that $(B^s_t, Q^{s,w}_t, W^{s,w}_t)$, $t\ge s$, is an $\hF^s$-adapted  process defined on $(\O, \cF, \hP_{sw})$, and it is Markovian.

\subsection{Optimal Dividend-Investment Problem with the Sparre Andersen Model}

In this paper we assume that the dynamics of surplus of an insurance company, denoted by $X=\{X_t\}_{t\ge0}$, in the absence of dividend payments and investment, is described by the following {\it Sparre Andersen model} on the given probability space
$(\O, \cF, \hP; \hF)$:
\bea
\label{X0}
X_{t}=x+pt-Q_t:=x+pt-\sum_{i=1}^{N_{t}}U_{i},\q t\in [0,T],
\eea
where  $x=X_{0}\geq 0$, $p>0$ is the constant premium rate, and $Q_t=\sum_{i=1}^{N_{t}}U_{i}$ is the (renewal) claim process.
%
%
We shall assume that the insurer is allowed to both invest its surplus in a financial market and will also pay dividends, and
will try to maximize the dividend received before the ruin time of the insurance company. To be more precise, we shall assume
that the financial market is described by the  standard Black-Scholes model. 
That is, the prices of the risk-free and risky assets satisfy the following SDE
\bea
 \label{market}
 \left\{\ba{lll}
  dS^0_{t}=rS^0_{t}dt, \\
 dS_{t}= \mu S_{t}dt+\sigma S_{t}dB_{t}, 
 \ea\right. \qq t\in[0,T],
 \eea
where $B=\{B_{t}\}_{t\ge0}$ is the given Brownian motion, $r$ is the interest rate, and $\mu>r$ is the 
appreciation rate of the stock. 

With the same spirit, in this paper we shall consider a portfolio with only one risky asset and one bank account
and define the
control process by $\pi=(\gamma_{t},L_{t})$, $t\ge0$, where $\g\in L^2_{\hF}([0,T])$ is a self-financing strategy, representing the proportion of the surplus invested in the stock at time $t$ (hence $\gamma_{t}\in [0,1]$, for all
$t\in[0,T]$), and $L\in L^2_{\hF}([0,T])$ is the cumulative dividends the company has paid out up to time $t$
(hence $L$ is increasing). Throughout this paper we will consider the the filtration $\hF=\hF^{(B,Q, W)}$, and we say that a control strategy $\pi=(\gamma_{t},L_{t})$ is
{\it admissible} if it is $\hF$-predictable with c\`adl\`ag paths, and square-integrable (i.e., $\hE\big[\int_0^T|\g_t|^2dt+|L_T|^2\big]<\infty$.)
and we denote the set of all admissible  strategies
restricted to $[s, t]\subseteq [0,T]$ by $\sU_{ad}[s,t]$. Furthermore, we shall often  use the notation $\sU_{ad}^{s,w}[s, T]$ to specify the probability space $(\O, \cF, \hP_{sw})$, and denote $\sU^{0,0}_{ad}[0,T]$ by $\sU_{ad}[0,T]=\sU_{ad}$ for simplicity.

By a standard argument using the self-financing property, one can easily show that, for any $\pi\in \sU_{ad}$ and any initial surplus ${x}$, the dynamics of the controlled risk process $X$ satisfies the following SDE:
\bea
\label{surplus}
dX_{t}^{\pi}=pdt+rX_{t}^{\pi}dt+(\mu-r)\gamma_{t}X^\pi_{t}dt+ \sigma \gamma_{t}X_{t}^{\pi}dB_{t}-dQ_t-dL_{t}, \q X^\pi_0=x, \q t\in [0,T].
\eea
We shall denote the solution to (\ref{surplus}) by $X_t=X^\pi_t=X_{t}^{\pi,x}$, whenever the specification of $(\pi, x)$  is necessary. Moreover, for any $\pi\in \sU_{ad}$, we denote $\tau^{\pi}=\t^{\pi,x}:={\rm inf}\{t\geq0; X_{t}^{\pi,x}< 0\}$ to be the ruin time of the insurance company.   We shall make use of the following {\it Standing Assumptions}:
\begin{assum}
\label{assump1}
(a) The interest rate $r$, the volatility $\si$, and the insurance premium $p$ are all positive constants,;

(b)  The distribution functions $F$ (of $T_i$'s) and $G$ (of $U^i$'s) are continuous on $[0,\infty)$. Furthermore, $F$ is absolutely continuous, with density function $f$ and intensity function $\l(t):=f(t)/\bar F(t)>0$, $t\in[0,T]$;

(c) The cumulative dividend process $L$ is absolutely continuous with respect to the Lebesgue measure. That is,
 there exists $a\in L^2_\hF([0,T]; \hR_+)$, such that $L_t=\int_0^t a_sds$, $t\ge 0$.
 We assume further that for some  constant $M\ge p>0$, it holds that $0\le a_t\le M$, $dt\times d\hP$-a.e.
\end{assum}

\begin{rem}
\label{remark1}
{\rm
 1)   Since in this paper we are focusing mainly on the value function and the dynamic programming approach, we can and shall assume that we are under the risk-neutral measure, that is, $\m=r$ in (\ref{surplus}). We should
 note that such a simplification does not change the technical nature of any of our discussion.

2) The Assumption \ref{assump1}-(c) is merely  technical, and it is not unusual, see for example, \cite{AT, Ger2, jea}. But this assumption will certainly exclude the possibility of having ``singular"
type of strategies, which could very well be the form of an optimal strategy in this kind of problem. However, since in
this paper our main focus is to deal with the difficulty caused by the renewal feature of the model, we are content with such an assumption.}
\qed
\end{rem}

We should note that the surplus dynamics (\ref{surplus}) with Assumption \ref{assump1}-(a) is in the simplest form.
More general dynamics with carefully posed assumptions is clearly possible, but not essential for the main results of this paper. In fact, as we can see later, even in this simple form the technical difficulties are already significant. We therefore
prefer not to pursue the generality of the surplus dynamics in the current paper so as not to disturb the already lengthy presentation.
In the rest of the paper we
shall consider, for given $s\in[0,T]$,  the following SDE (recall (\ref{surplus}) and Remark \ref{remark1} on the filtered probability space $(\O, \cF, \hP_{sw}; \hF^s)$):
for $(\g, a)\in \sU_{ad}^{s,w}[s, T]$,
%
\bea
\label{Xsxw}
\left\{\ba{lll}
\dis X_{t}^{\pi}=x+p(t-s)+r\int_{s}^{t}X_{u}^{\pi}du+\sigma \int_{s}^{t}\gamma_{u}X_{u}^{\pi}dB_{u}-Q^{s,w}_t-\int_s^t
a_udu, \\
W_t :=w+(t-s)-(\si_{N_t}-\si_{N_s}),
\ea\right.
\q t\in [s,T].
\eea
We denote the solution by $(X^\pi,W)=(X^{\pi,s,x,w}, W^{s,w})$, to emphasize its dependence on the initial state $(s,x,w)$.

We now describe our optimization problem. Given an admissible strategy $\pi\in \sU_{ad}^{s,w}[s, T]$, we define the {\it cost functional}, for the given initial data  $(s,x,w)$ and
%
the  state  dynamics (\ref{Xsxw}), as
\bea
\label{cost1}
J(s,x,w;\pi);=\hE_{sw}\Big\{\int_{s}^{\tau^{\pi}_s\wedge T}e^{-c (t-s)}dL_{t}\Big|X^\pi_{s}=x\Big\}:=\hE_{sxw}\Big\{\int_{s}^{\tau^{\pi}_s\wedge T}e^{-c (t-s)}dL_{t}\Big\}.
\eea
Here $c>0$ is the {\it discounting factor} (or {\it force of interest}), and $\t^\pi_s=\t^{\pi,x,w}_s:=\inf\{t >s: X^{\pi,s,x,w}_{t}<0\}$ is
the ruin time of the insurance company. That is, $J(s,x, w;\pi)$ is the expected total discounted amount of dividend received until the ruin. Our objective is to find the optimal strategy $\pi^*\in\sU_{ad}[s, T]$ such that
 \bea
 \label{V1}
 V(s,x,w):= \sup_{\pi\in \sU_{ad}[s, T]} J(s,x,w; {\pi}).
 \eea

%
%

%


%
 We note that the value function should be defined for $(s,x,w)\in D$ where $D=\{(s,x,w): 0\leq s \leq T, x\geq0, 0\leq w \leq s\}$. We make the convention that  $V(s,x,w)=0$, for $(s,x,w)\notin D$. We shall frequently carry out our discussion on the following two sets:
 \bea
 \label{D}
 \sD&:=&\mbox{int}D=\{(s,x,w)\in D: 0<s<T, \, x>0,\, 0<w<s\}; \\
  \sD^*&:=&\{(s,x,w)\in D: 0\leq s<T,\, x\geq0, \, 0\leq w\leq s\}.\nonumber
   \eea
  We note that $\sD\subset \sD^*\subseteq \bar\sD=D$, the closure of $\sD$, and $\sD^*$ does not include boundary at the terminal time $s=T$.

To end this section we list two technical lemmas that will be useful in our discussion. The proofs of these lemmas are very similar to
the Brownian motion case (cf. e.g., \cite[Chapter 3]{yong-zhou}), along the lines of {\it Monotone Class Theorem} and {\it Regular
Conditional Probability Distribution} (RCPD), we therefore
omit them.  Let us denote $D^m_T:=D([0,T];\hR^m)$, the space of all $\hR^m$-valued c\`adl\`ag functions on $[0,T]$, endowed with the sup-norm, and $\sB^m_T:=\sB(D^m_T)$, the topological Borel field on $D^m_T$. Let $D_t^m:=\{\zeta_{\cdot\wedge t}| \zeta\in D^m_T\}$,
$\sB_t^m:= \sB(D_t^m)$, $ t\in [0,T]$, and $\sB_{t^+}^m:=\cap_{s>t}\sB_s^m$, $t\in[0,T]$.  For a generic Euclidean space
  $\hX$,
we denote $\sA_T^m(\hX)$ to be the set of all $\{\sB_{t^+}^m\}_{t\geq 0}$-progressively measurable process $\eta: [0,T]\times D^m_T\rightarrow \hX$. That is, for any $\f \in \sA^m_T(\hX)$, it holds that $\f(t, \eta)=\f(t, \eta_{\cd\wedge t})$, for $t\in[0,T]$ and
$\eta\in D^m_T$. As usual, we denote $\sA^m_T=\sA^m_T(\hR)$ for simplicity.
\begin{lem}
\label{A1}
Let $(\Omega,{\cal F}, \hP)$ be a complete probability space, and  $\zeta: \Omega\rightarrow D_T^m$ be a $D^m$-valued
process. Let ${\cal F}_t^{\zeta}=\sigma\{\zeta(s):0\leq s\leq t\}$. Then $\phi:[0,T]\times \O\mapsto \hX$ is
$\{{\cal F}_t^{\zeta}\}_{t\geq 0}$-adapted if and only if there exists an $\eta\in \sA_T^m(\hX)$ such that $\phi(t,\omega)=\eta(t,\zeta_{\cdot\wedge t}(\omega))$, $ \hP$-a.s.-$\o\in \O$, $ \forall t\in [0,T]$.
\qed
\end{lem}

\begin{lem}
	\label{PDPP}
Let  $(s,x,w)\in D$ and $\pi=(\g, a)\in {\sU}_{ad}[s,T]  $. Then for any stopping time $\t\in[s, \t^\pi]$, $\hP$-a.s., and any
${\cal F}_{\t}$-measuable random variable $(\xi,\eta)$ taking values in $[0,\infty)\times[0,T]$, it holds that
	\bea
	J(\t,\xi(\omega),\eta(\omega);\pi)=\hE\Big\{\int_{\t}^{\t^\pi}e^{-c(t-\t)}a_tdt\Big|\cF_{\t}\Big\}(\omega), ~\mbox{ for $\hP$-a.s. $\o\in\O$.}
	\eea
	\qed
	\end{lem}

 \section{Basic Properties of the Value Function}
\setcounter{equation}{0}

In this section, we present some results that characterize the regularity of the value function $V(s,x,w)$. We should note
that the presence of the renewal process, and consequently the extra component $W=\{W_t\}_{t\ge0}$, changes the nature
of the dynamics significantly. In fact, even in this simple setting, many well-known properties of the value function becomes
either invalid, or much less obvious.

We begin by making some simple but important observations, which will be used throughout the paper. First, we note that
 in the absence of claims (or in between the jumps of $N$), for a given $\pi=(\g, a)\in \sU_{ad}[s,T]$, the
dynamics of the surplus follows a non-homogeneous linear SDE (\ref{Xsxw}) with $Q^{s,w}\equiv0$, and
has the explicit form (cf. \cite[p.361]{KS1}):
\bea
\label{Xsol}
X^\pi_t=Z^s_t\Big[X^\pi_s+\int_s^t [Z^s_u]^{-1}(p-a_u)du\Big], \q t\in[s, T],
\eea
where $Z_t^s :=\exp\Big\{r(t-s)+\si\int_s^t\g_udB_u-\frac{\si^2}2\int_s^t|\g_u|^2du\Big\}$. From (\ref{X0}) and (\ref{Xsol}) it is
clear that in the absence of claims, the surplus $X_t<0$ would never happen if one does not over pay the dividend  whenever $X_t=0$. For example, if we consider only those  $\pi\in \sU_{ad}$ such that
%
$(p-a_t){\bf 1}_{\{X^\pi_t=0\}}\ge 0$, $\hP$-a.s., then we have $dX^\pi_t\ge 0$, whenever $X^\pi_t=0$, which implies that $X^\pi_t\ge 0$ holds
for all $t\ge 0$.
Such an assumption, however, would cause some unnecessary complications on the well-posedness of the SDE (\ref{X0}). We
shall argue slightly differently.

Since it is intuitively clear that
the dividends should only be paid when reserve is positive, we suspect that any $\pi\in \sU_{ad}$ such that $\t^\pi$ occurs
in between claim times (i.e., caused by overpaying dividends) can never be optimal.
The following result justifies this point.
\begin{lem}
\label{ruin}
Suppose that $\pi\in\sU_{ad}$ is such that $\hP\{\si_i\wedge T<\t^\pi<\si_{i+1}\wedge T\}>0$, for some $i\in\hN$,
where $\si_i$'s are the jump
times of $N$, then there exists $\tilde\pi\in\sU_{ad}$
such that $\hP\{\t^{\tilde \pi}\in \bigcup_{i=1}^\infty \si_i\}=1$, and $J(s,x,w;\tilde \pi)>J(s,x,w;\pi)$.
\end{lem}

{\it Proof.} Without loss of generality we assume $s=w=0$. We first note from (\ref{Xsol}) that on the set $\{\si_i\wedge T<\t^\pi<\si_{i+1}\wedge T\}$, one must have $X^\pi_{\t^\pi-}=X^\pi_{\t^\pi}=0$, and for
some $\d>0$, $a_t>p$ for $t\in[\t^\pi, \t^\pi+\d]$. Now define $\tilde\pi_t:= \pi_t{\bf 1}_{\{t<\t^\pi\}}+(0,p){\bf 1}_{\{t\ge \t^\pi\}}$, and
denote $\tilde X:=X^{\tilde\pi}$. Then clearly, $\tilde X_t=X^\pi_t$ for all $t\in [0,\t^\pi]$, $\hP$-a.s, and $d\tilde X_{\t^\pi}=(p-\tilde a_{\t^\pi})dt=0$. Consequently $\tilde X_t\equiv 0$ for $t\in [\t^\pi, \si_{i+1}\wedge T)$ and $\tilde X_{\si_{i+1}}<0$ on $\{\si_{i+1}<T\}$. In other words,
$\t^{\tilde \pi}=\si_{i+1}$, and thus
\beaa
J(0,x,0;\tilde \pi)&=&\hE\Big[\int_0^{\t^{\tilde \pi}\wedge T}e^{-ct}a_tdt\Big]\\
&\geq &J(0,x,0;\pi)+\hE\Big[\int_{\t^\pi}^{\si_{i+1}\wedge T} p
e^{-ct}dt: \si_{i}\wedge T<\tau^{\pi}<\si_{i+1}\wedge T\Big]> J(0,x,0;\pi),
\eeaa
since  $\hP\{\si_i\wedge T<\t^\pi<\si_{i+1}\wedge T\}>0$, proving the lemma.
\qed

We remark that Lemma \ref{ruin} amounts to saying that for an optimal policy it is necessary that
ruin only occurs at the arrival of a claim. Thus, in the sequel we shall consider a slightly fine-tuned set of
{\it admissible stratetigies}:
\bea
\label{admL1}
\tilde \sU_{ad}:=\Big\{\pi=(\g, a)\in\sU_{ad}: \D X^\pi_{\t^\pi}{\bf 1}_{\{ \t^\pi<T\}}<0,    \hP\mbox{-a.s.}\Big\}.
\eea
The set $\tilde\sU_{ad}[s,T]$ is defined similarly for $s\in[0,T]$, and we shall often drop the ``$~\tilde{}~$" for simplicity.

%
%

We now  list some generic properties of the value function.
\begin{prop}
\label{vb}
Assume that the Assumption \ref{assump1} is in force. Then, the value function $V$ enjoys the following properties:

  (i) $V(s,x,w)$ is increasing with respect to $x$;

   (ii) $V(s,x,w)\leq \frac{M}{c}(1-e^{-c(T-s)})$ for any $(s,x,w)\in D$, where $M>0$ is the
   constant given in Assumption \ref{assump1}; and

   (iii) $\lim_{x\rightarrow\infty}V(s,x,w)=\frac{M}{c}[1-e^{-c(T-s)}]$, for $0\leq s \leq T$, $ 0\leq w \leq s $.

\end{prop}
{\it Proof.} (i) is obvious, given the form of the solution (\ref{Xsol}); and (ii) follows from the simple estimate:
$V(s,x,w)\leq \int_{s}^{T}e^{-c(t-s)}Mdt=\frac{c}{M}[1-e^{-c(T-s)}]$.

To see (iii), we consider a simple strategy:
$\pi^0:=(\gamma,a)\equiv (0,M)$.
Then we can write
\bea
\label{Xpi0}
X^{\pi^0,x,w}_t=e^{r(t-s)}x+\frac{p-M}{r}\Big(1-e^{-r(t-s)}\Big)-\int_s^t e^{r(t-u)}dQ^{s,w}_u, \q t\in[s,T],
\eea
and it is obvious that
$\lim_{x\to\infty} \t^{\pi^0,x,w}_s\wedge T=T$, $\hP$-a.s. Thus we have
\beaa
V(s,x,w)&\ge& J(s,x,w;{\pi^0})=\hE\Big[\int_{s}^{\tau^{\pi^0,x,w}_s\wedge T}e^{-c(t-s)}Mdt\Big]= \frac{M}{c}\hE\Big[1-e^{-c(\t^{\pi^0,x,w}_s\wedge T-s)}\Big].
\eeaa
By the Bounded Convergence Theorem we have
$\lim_{x\to\infty}V(s,x,w)\ge \frac{M}{c}(1-e^{-c(T-s)})$. This, combined with (ii), leads to (iii).
\qed


In the rest of this subsection we study the continuity of the value function $V(s,x,w)$ on the temporal variable $s$, for fixed initial state $(x,w)$.  We have the following result.
\begin{prop}
\label{conts}
Assume Assumption \ref{assump1}. Then, $\forall (s,x,w), (s+h, x,w)\in D$, $h>0$, it holds

(a) $V(s+h,x,w)-V(s,x,w)\leq 0$;

(b) $V(s,x,w)-V(s+h,x,w)\leq Mh$, where $M>0$ is the constant in Assumption \ref{assump1}.
\end{prop}
{\it Proof.}  We note that the main difficulty here is that, for the given $(s,x,w)$, the claim process $Q^{s,w}_t=\sum_{i=1}^{N^{s,w}_t}U_i$, $t\ge0$, and the ``clock" process
 $W^{s,w}=\{W^{s,w}_t\}_{t\ge0}$ cannot be controlled, thus it is not possible to keep the process $W$ ``frozen" at the initial
 state $w$
 during the time interval $[s,s+h]$ by any control strategy. We shall try to get around this by using a
``time shift" to move the initial time to $s=0$.

More precisely, for any $\pi\in {\sU}_{ad}^{s,w}[s,T]$, we define $\bar{\pi}^s_t=(\bar{\gamma}_{t}, \bar{a}_{t}):=(\gamma_{s+t}, a_{s+t})$, $t\in [0,T-s]$. Then $\bar\pi$ is adapted to the filtration $\bar\hF^s:=\{\cF_{s+t}\}_{t\ge0}$.
consider the optimization problem on the new probability set-up $(\O, \cF, \hP_{sw}, \bar\hF^{s}; \bar B^s, \bar Q^{s,w}, \bar W^{s,w})$, where $(\bar B^s_t, \bar {Q}^{s,w}_{t},\bar W^{s,w}_{t})=(B^s_{s+t},   {Q}^{s,w}_{s+t},W^{s,w}_{s+t})$, $t\ge0$.
Let us
denote the corresponding admissible control set by  $\bar{\sU}^{s,w}_{ad}[0,T-s]$, to emphasize the obvious dependence on
the initial state $(s,w)$. Then
$\bar\pi^{s}\in \bar {\sU}^{s,w}_{ad}[0,T-s]$,  and the corresponding surplus process, denoted by $\bar{X}^{\bar{\pi}^{s}}$, should satisfy the SDE:
 \bea
 \label{bXsxw}
\bar{X}_{t}^{\bar{\pi}^{s}}=x+pt+r\int_{0}^{t}\bar{X}_{u}^{\bar{\pi}^{s}}du
+\sigma \int_{0}^{t}\bar{\gamma}_{u}\bar{X}_{u}^{\bar{\pi}^{s}}d\bar{B}^s_u-\bar Q^{s,w}_{t}-\int_{0}^{t}\bar{a}_{u}du, \q t\ge 0.
\eea
Since the SDE is obviously pathwisely unique, whence unique in law, we see that the laws of  $\{\bar X^{\bar{\pi}^{s}}_t\}_{t\ge0}$  and that of $\{X_{s+t}^\pi\}_{t\ge0}$ (which satisfies (\ref{Xsxw})), under $\hP_{sw}$, are identical. In other words, if we specify
the time duration in the cost functional, then we should have
\bea
\label{Jsw}
\left\{\ba{lll}
\dis J_{s,T}(s,x,w;\pi):=\hE_{sw}\Big[\int_{s}^{\tau^{{\pi}}\wedge T}e^{-c(t-s)}a_{t}dt|X^\pi_{s}=x\Big]\ms\\
\qq\qq\dis =\hE_{sw}\Big[\int_{0}^{\tau^{\bar{\pi}^s}\wedge (T-s)}e^{-ct}\bar{a}_{t}dt|\bar{X}^{\bar\pi}_{0}=x\Big]=:\bar J_{0,T-s}(0,x,w;\bar\pi^s),\ms\\
\dis V(s,x,w)=\sup_{\bar\pi\in \bar\sU^{s,w}_{ad}[0,T-s] }\bar J_{0,T-s}(0,x,w;\bar\pi).
\ea\right.
\eea
Similarly, for any $\pi\in {\sU}_{ad}[s+h,T]$, we can find $\hat\pi^{s+h}\in {\sU}^{s+h,w}_{ad}[0,T-s-h] $ such that
\bea
\label{Jshw}
\left\{\ba{lll}
J_{s+h, T}(s+h,x,w,\pi)=\hat J_{0,T-s-h}(0,x,w,\hat\pi^{s+h}), \ms\\
\dis V(s+h,x,w)=\sup_{\hat\pi\in \bar \sU^{s+h,w}_{ad}[0,T-s-h] }\hat J_{0,T-s-h}(0,x,w;\hat\pi).
\ea\right.
\eea
Now, for
the given $\hat{\pi}\in\bar{\sU}^{s+h,w}_{ad}[0,T-s-h]$ we apply Lemma \ref{A1} to find $\eta\in\sA^3_{T-s-h}(\hR^2)$,
such that $\hat{\pi}_{t}=\eta(t,\bar B^{s+h}_{\cdot\wedge t}, \bar Q^{s+h,w}_{\cdot\wedge t},\bar W^{s+h,w}_{\cdot\wedge t})$, $ t\in[0,  T-s-h]$.
We now define
$$ \tilde{\pi}^h_t:=\eta(t, \bar B^s_{\cdot\wedge t\wedge (T-s-h)}, \bar Q^{s,w}_{\cdot\wedge t\wedge(T-s-h)},\bar W^{s,w}_{\cdot\wedge t\wedge(T-s-h)}), \qq t\in [0,T-s].
$$
Then, $\tilde \pi^h\in \bar\sU^{s,w}_{ad}[0,T-s]$. Furthermore, since the law of $(\bar B^{s+h}_t, \bar Q^{s+h,w}_t,\bar W^{s+h,w}_t)$, $t\in[0,T-s-h]$, under $\hP_{(s+h)w}$,  and that of $(\bar B^{s}_t, \bar Q^{s,w}_t,\bar W^{s,w}_t)$, $t\in[0,T-s-h]$, under $\hP_{sw}$, are identical, by the pathwise uniqueness (whence uniqueness in law) of the solutions to SDE (\ref{Xsxw}), the processes $\{(X^{\tilde\pi^h}_t,\bar W^{s, w}_t, \tilde\pi^h_t)\}_{t\in[0,T-s-h]}$ and
$\{(X^{\hat\pi}_t, \bar W^{s+h, w}_t, \hat\pi_t)\}_{t\in[0,T-s-h]}$ are identical in law. Thus
\beaa
 J_{s+h, T}(s+h,x,w,\pi)&=&\hat J_{0,T-s-h}(0,x,w;\hat\pi)
=\hE_{0xw}\Big[\int_0^{\t^{\hat\pi\wedge (T-s-h)}}e^{-ct}\hat a_tdt|\Big]\\
&\le& \hE_{0xw}\Big[\int_0^{\t^{\bar\pi\wedge (T-s)}}e^{-ct}\bar a_tdt  \Big]
= \bar J_{0,T-s}(0,x,w;\bar \pi)
 \le V(s,x,w).
 \eeaa
 Since $\pi\in \sU_{ad}[s+h,T]$ is arbitrary, we obtain  $ V(s+h,x,w)\le V(s,x,w)$,  proving (a).

To prove (b),  let $\pi\in\sU_{ad}[s,T]$. For any $h\in (0,T-s)$, we 
define $\pi^h_t:=\pi_{t-h}$ for $t\in [s+h,T]$. Then clearly, $\pi^h\in\sU_{ad}^{s+h, w}[s+h,T]$. 
Furthermore, we have
\bea
\label{partb}
&&J(s,x,w;\pi)-J(s+h,x,w;\pi^h)\nonumber\\
 &=&\hE_{sxw}\Big[\int_s^{\t^{\pi}}e^{-c(t-s)} a_tdt:\t^{\pi}\leq T-h\Big]+\hE_{sxw}\Big[\int_s^{\t^{\pi}\wedge T}e^{-c(t-s)} a_tdt:\t^{\pi}>T-h\Big]\nonumber\\
 &&-\hE_{(s+h)xw}\Big[\int_{s+h}^{\t^{\pi}}e^{-c(t-s-h)} a_{t-h}dt :\t^{\pi^h}\leq T\Big]-\hE_{(s+h)xw}\Big[\int_{s+h}^{T}e^{-c(t-s-h)} a_{t-h}dt :\t^{\pi^h}>T\Big].\nonumber
 \eea
By definition of the strategy $\pi^h$, it is easy to check that 
\beaa
\hE_{sxw}\Big[\int_s^{\t^{\pi}}e^{-c(t-s)} a_tdt:\t^{\pi}\leq T-h\Big]&=&\hE_{(s+h)xw}\Big[\int_{s+h}^{\t^{\pi}}e^{-c(t-s-h)} a_{t-h}dt :\t^{\pi^h}\leq T\Big]\\
\hE_{sxw}\Big[\int_s^{T}e^{-c(t-s)} a_tdt:\t^{\pi}>T-h\Big]&=&\hE_{(s+h)xw}\Big[\int_{s+h}^{T}e^{-c(t-s-h)} a_{t-h}dt :\t^{\pi^h}>T\Big],
   \eeaa
 we deduce from (\ref{partb}) that
\bea
\label{partb1}
J(s,x,w;\pi)-J(s+h,x,w;\pi^h)\le \hE_{sxw}\Big[\int_{T-h}^{T}e^{-c(t-s)} a_tdt
\Big]\leq Mh.
\eea
Consequently, we have $J(s, x,w;\pi)\le Mh+ V(s+h, x,w)$. 
Since $\pi\in \sU_{ad}^{s,w}[s,T]$ is arbitrary, we obtain
 (b), proving the proposition.
\qed

We complete this section with an estimate that is quite useful in our discussion. First note that (\ref{Xsol}) implies that
in the absence of claims, the surplus without investment and dividend (i.e., $\pi\equiv (0,0)$) is $X^{0,s, x}_t= e^{r(t-s)}[x+\frac{p}{r}(1-e^{-r (t-s)})]$.
\begin{prop}
\label{1sLip}
Let $(s,x,w)\in D$. Then, for any $h>0$ such that $(s+h, X^{0,s,x}_{s+h}, w+h)\in D$, it holds that
\bea
\label{Vest}
V(s+h,X^{0,s,x}_{s+h},w+h)\leq e^{ch+\int_{w}^{w+h}\frac{f(u)}{\bar{F}(u)}du}V(s,x,w).
\eea
\end{prop}

{\it Proof.}   For any $\varepsilon>0$, we choose $\pi^{h,\e}\in {\sU}^{s+h, w+h}_{ad}[s+h,T]$  such that
$$J(s+h,X^{0,s,x}_{s+h},w+h;\pi^{h, \e})\geq V(s+h,X^{0,s,x}_{s+h},w+h)-\varepsilon.$$

%
Now define a new strategy: $\bar \pi^{h,\e}_t=\pi^{h,\e}_t{\bf 1}_{\{T_1^{s,w}> h\}}
 {\bf 1}_{[s+h,T]}(t)$, $t\in[s, T]$, where
$T^{s,w}_1$ is the first jump time of the delayed renewal process $N^{s,w}$.
Then, clearly, $\bar{\pi}^{h, \e}\in{\sU}^{s,w}_{ad}[s,T]$,
 and
 $X_{s+h}^{\bar{\pi}^h}=X^{0,x}_{s+h}$ on the set $\{T_1^{s,w}>h\}\in\cF_{s+h}$.
Thus, using (\ref{PTsw}) we have
\begin{eqnarray*}
V(s,x,w)&\geq& J(s,x,w;\bar{\pi}^{h,\e})=\hE_{sxw}\Big[\int_{s+h}^{\bar\t^{h,\e}\wedge T} e^{-c(t-s)}a^{h,\e}_tdt \cd{\bf 1}_{\{
T_1^{s,w}>h\}}\Big]\\
&= & e^{-ch}J(s+h,X^{0,s,x}_{s+h},w+h;\pi^{h,\e})\hP_{sxw}\{T_1^{s,w}>h\}\\
&\ge& [V(s+h,X^{0,s,x}_{s+h},w+h)-\e]e^{-ch-\int_{w}^{w+h}\frac{f(u)}{\bar{F}(u)}du}.
\end{eqnarray*}
Letting $\e\to0$ we obtain the result.
\qed

We note that a direct consequence of (\ref{Vest}) is the following inequality:
\bea
\label{Vest1}
V(s+h,X^{0,s,x}_{s+h},w+h)-V(s,x,w)\leq [e^{ch+\int_{w}^{w+h}\frac{f(u)}{\bar{F}(u)}du}-1]V(s,x,w).
\eea
This gives a kind of one-sided continuity of the value function, although it is a far cry from a true
joint continuity which we will study in the next sections.

\section{Continuity of the value function on $x$}
\setcounter{equation}{0}

In this section we investigate the continuity of value function on initial surplus $x$. As in all ``exit-type" problem, the main
subtle point here is that the ruin time $\t^\pi$, which obviously depend on the initial state $x$, is generally not continuous
in $x$. We shall  borrow the idea of {\it penalty method} (see, e.g., \cite{FS}), which we now describe.

To begin with, we recall the domain  $D=\{(s,x,w)\neg:\neg 0\le s\le T, x\ge 0, 0\le w\le s\}$. Let $d(x,w):=(-x)\vee 0$ for
$(x,w)\in \hR\times[0,T]$, and for $\pi\in\sU^{s,w}_{ad}[s,T]$ we define a {\it penalty function} by
\bea
\label{beta}
\beta(t,\e)=\beta^{\pi,s,x,w}(t,\e)=\exp\Big\{-\frac{1}{\e}\int_s^t d(X_r^{\pi,s,x,w},W^{s,w}_r)dr\Big\}, \qq t\ge 0.
\eea
Then clearly $\beta(t, \e)=1$ for $t\le \t^\pi_s$. Thus we have
\bea
\label{VegeV}
V^{\e}(s,x,w)&=&\sup_{\pi\in \sU_{ad}[s,T]}J^{\e}(s,x,w;\pi):=\sup_{\pi\in \sU_{ad}[s,T]}\hE\Big[\int_s^T \beta^{\pi,s,x,w}(t,\e)e^{-c(t-s)}a_tdt\Big]\\
&=&\sup_{\pi\in \sU_{ad}[s,T]}\hE\Big[\int_s^{\t^\pi_s} e^{-c(t-s)}a_tdt+\int_{\t^{\pi}_s}^T \beta^{\pi,s,x,w}(t,\e)e^{-c(t-s)}a_tdt\Big]\ge V(s,x,w).
\nonumber
\eea
We have the following lemma.
\begin{lem}
\label{Vecont}
Let $K\subset D$ be any compact set. Then the mapping $x\mapsto V^{\e}(s,x,w)$ is continuous, uniformly for  $(s,x,w)\in K$.
\end{lem}

{\it Proof.}
For $\pi\in\sU^{s,w}_{ad}[s,T]$, and $x_1, x_2\in[0,\infty)$ we
have
\bea
\label{estVe}
&& \hE\left|\beta^{\pi,x_1}(t,\e)-\beta^{\pi,x_2}(t,\e)\right| \nonumber\\
&= &\hE \left|\exp\Big\{-\frac{1}{\e}\int_s^t d(X^{\pi,x_1}_r,W_r)dr\Big\}-\exp\Big\{-\frac{1}{\e}\int_s^t d(X^{\pi,x_2}_r,W_r)dr\Big\}\right| \notag\\
 &\le&\frac{1}{\e} \hE\left|\int_s^t d(X^{\pi,x_1}_r,W_r)-d(X^{\pi,x_2}_r,W_r)dr\right|\le \frac{1}{\e}\int_s^t \hE |(X^{\pi,x_1}_r-X^{\pi,x_2}_r| dr\\
 &\le &\sqrt{T}\frac{1}{\e} (\int_s^t \hE\left|X^{\pi,x_1}_r-X^{\pi,x_2}_r\right|^2dr)^{\frac{1}{2}}
\le \frac{T}{\e}|x_1-x_2|. \nonumber
\eea
%
In the above, the last inequality is due to a standard estimate of the SDE (\ref{X0}).
Thus, by some standard argument, we conclude that $V^{\e}$ is continuous in $x$. Since $K$ is compact,
the continuity is uniform for $(s,x,w)\in K$.
\qed

We should point out that, the estimate (\ref{estVe}) indicates that the continuity of $V^\e$ (in $x$), while uniformly on compacta,
is not uniform in $\e$(!). Therefore, we are to argue that, as $\e\to0$, $V^\e\to V$ on any compact set $K\subset D$, and the
convergence is uniform in all $(s,x,w)\in K$, which would in particular imply  that $V$ is continuous  on $D$. In other words, we are aiming
at the following main result of this section.
\begin{thm}
\label{continuityx}
For any compact set $K\subset D$, the mapping $x\mapsto V(s,x,w)$ is continuous, uniformly for $(s,x,w)\in K$. In particular,
the value function $V$ is continuous in $x$, for $x\in[0,\infty)$.
\end{thm}

To prove Theorem \ref{continuityx}, we shall introduce an intermediate problem. For each $\th>0$, we denote $D_{\th}
:=\{(s, x,w): s\in [0,T], x\in(-\th,\infty),w\in[0,s]\}$.
Clearly $D_\th\subset D_{\th'}$ for $\th<\th'$, and $\bigcap_{\th>0}D_\th=D$.  For  $(s,x,w)\in K$ and
$\pi\in\sU_{ad}[s,T]$, we denote
$\t^{\pi, \th}_s=\t^{\pi,\th}_{s,x,w}$ (resp. $\t^{\pi,0}_s$) to be the exit time of the process $(t,X^{\pi, s,x,w}_t, W^{s,w}_t)$ from $D_\th$ (resp. $D$) before $T$.
%
For notational simplicity we shall write  $(X^\pi,W):=(X^{\pi,s,x,w},W^{s,w})$, $\t:=\t^{\pi,0}_s$, and $\t^\th:=\t^{\pi, \th}_{s}$, when
the context is clear.  It is worth noting that the  function $\beta(t,\e)$ satisfies an SDE:
\beaa
\beta(t,\e)=1-\frac{1}{\e}\int_s^td(X^{\pi}_r,W_r)\beta(r,\e)dr, \qq t\in[s,T].
\eeaa
Thus, together with the underlying process $(X^{\pi},W)$, we see that the optimization problem in
(\ref{VegeV}) is
a standard stochastic control problem with jumps and fixed terminal time $T$,  therefore
the standard Dynamic Programming Principle (DPP) holds for $V^\e$. To be more precise, for any
stopping time $\hat\t\in[s,T]$, it holds that
\bea
\label{DPPVe}
 V^{\e}(s,x,w)=\sup_{\pi\in \sU_{ad}[s,T]}\hE_{sxw}\Big\{\int_s^{\hat\t}\beta(t,\e)e^{-c(t-s)}a_tdt+e^{-(\hat\tau-s)}\beta(\hat\tau,
 \e) V^{\e}(\hat\t,X^\pi_{\hat\t},W_{\hat\t})\Big\}.
 \eea

We are now ready to prove Theorem \ref{continuityx}.

\ms
\no[{\it Proof of Theorem \ref{continuityx}.}] We first note that, for any $(s,x,w)\in K$ and $\pi\in\sU_{ad}[s,T]$,  by DPP (\ref{DPPVe}) and the fact (\ref{VegeV}) we have
 \bea
 \label{VleVe}
&& V(s,x,w)\le V^{\e}(s,x,w)\notag\\
&=&\sup_{\pi\in \sU_{ad}[s,T]}\hE_{sxw}\{\int_s^{\t^\th}\beta(t,\epsilon)e^{-c(t-s)} a_{t}dt+e^{-(\t^\th-s)}\beta(\t^\th, \e)V^{\e}(\t^\th,X_{\t^\th},W_{\t^\th})\}\nonumber\\
&=&\sup_{\pi\in \sU_{ad}[s,T]}\hE_{sxw}\Big\{\int_s^{\t}e^{-c(t-s)} a_{t}dt+\int_{\t}^{\t^\th}\beta(t,\e)e^{-c(t-s)}a_{t}dt\\
&&\qq\qq\q+e^{-(\t^\th-s)}\beta(\t^\th, \e)V^{\e}(\t^\th,X^\pi_{\t^\th},W_{\t^\th})\}\notag\\
&\le& V(s,x,w)+C\sup_{\pi\in \sU_{ad}[s,T]}\hE_{sxw}(\t^\th-\t)+h^\th(\e),\notag
 \eea
where $h^\th(\e):=\hE_{sxw}[V^{\e}(\t^\th,X^\pi_{\t^\th},W_{\t^\th})]$, and $C>0$ is a generic constant depending only on the constants in
Assumption \ref{assump1} and $T$. We first argue that $\sup_{\pi\in \sU_{ad}[s,T]}\hE_{sxw}|\tau-\tau^\th|\to0$, as $\th\to 0$,  uniformly in $(s,x,w)\in K$.

 To see this, first note that
$ \sup_{\pi\in \sU^{s,w}_{ad}[s,T]}\hE_{sxw}|\tau-\tau^\th|
 \le\sup_{\pi\in \sU^{s,w}_{ad}[s,T]}T \hP\{\tau\ne\tau^\th\}$, here and in what follows  $\hP:=\hP_{sxw}$, if there is no danger of
 confusion. On the other hand, recall that $\t$ must happen at a claim arrival
 time on $\{\t\neq \t^\th\}$,
 and $\D X^\pi_t=\D Q_t^{s,w}$, it is easy to check that
 \beaa
 \label{Ptnoth}
 \hP\{\t\ne\t^\th\}
 &=&\hP\{\Delta X^\pi_{\t}\in (X^\pi_{\t-},X^\pi_{\t-}+\th)\}\\
 &=&\int_0^\infty \hP\Big\{\Delta Q_\t^{s,w}\in(y,y+\th)\Big|X^\pi_{\t-}=y\Big\}F_{X^\pi_{\t-}}(dy)=\int_0^\infty [G(y+\th)-G(y)]F_{X^\pi_{\t-}}(dy),
 \nonumber
  \eeaa
where $G$ is the common distribution function of the claim sizes $U^i$'s. Since $G$ is uniformly continuous on $[0,\infty)$, thanks to Assumption \ref{assump1}-(b), for any $\eta>0$ we can find $\th_0>0$, depending only on $\eta$,
such that $|G(y+\th_0)-G(y)| <\frac{\eta}{2T}$, for all $y\in[0,\infty)$,
\bea
\label{Etht}
\sup_{\pi\in \sU_{ad}[s,T]} \hE_{sxw}\{|\t^\th-\t|\}\le \sup_{\pi\in \sU_{ad}[s,T]}T\int_0^\infty |G(y+\th_0)-G(y)|F_{X^\pi_{\t-}}(dy)<\frac{\eta}2.
\eea

Plugging (\ref{Etht}) into (\ref{VleVe}) we obtain that
\bea
\label{VleVe1}
V(s,x,w)&\le& V^{\e}(s,x,w)\le V(s,x,w) +\frac{\eta}2+h^{\th_0}(\e).
\eea
We claim that  $ \lim_{\e\to0}h^{\th_0}(\e)=0$,
 and that the limit is uniform in $(s,x,w)\in K$.
 To this end,
 we define, for the given $\pi\in \sU_{ad}[s, T]$, and $\th=\th_0$,
 \bea
 \label{tM}
  \bar\t_\th:=\inf\{t>\t^\th,d(X_t^{\pi},W_t)<\th/2\}\wedge T;  \q
\bar \t_\th^{c}:=\inf\{t>\t^\th,d(X_t^{\pi,\th,c}, W_t)<\th/2\}\wedge T,
 \eea
 where $X^{\pi,\th,c}$ is the continuous part of $X^{\pi}$, for $t\ge\t^\th$, given $X^{\pi, \th,c}_{\t^\th}=X^{\pi}_{\t^\th}$.
Since $X^{\pi}$ only has negative jumps, we have  $\Delta X_t^{\pi}\le 0$, $\forall t\in[0,T]$.
Thus $\bar\t_\th^{c}\le\bar\t_\th$ and $d(X_t^{\pi, \th,c},W_t)\le d(X_t^{\pi},W_t)$, for all $t\in[s,T]$, $\hP$-a.s.  Furthermore, we note that
 $ d(X^{\pi,\th, c}_t, W_t)\ge \frac\th{2}$ for $t\in [\t^\th, \bar \t_\th^{c}]$, $\hP$-a.s.

 Now, denoting  $\hE_{\t^\th}[\,\cd\,]:=\hE[\,\cd\,| {\cal F}_{\tau_{\theta}}]$ and $X^c=X^{\pi,\th,c}$ we have, $\hP$-almost surely,
 \bea
 \label{Jep}
 &&J^{\e}(\t^\th,X_{\t^\th}^{\pi},W_{\t^\th};\pi)=\hE_{\tau^\th }\Big[\int_{\t^\th}^Te^{-\frac1\e\int_{\t^\th}^t d(X^\pi_r, W_r)dr}e^{-c(t-\t^\th)}a_tdt\Big]\nonumber\\
 &\le& \hE_{\tau^\th }\Big[\int_{\t^\th}^Te^{-\frac1\e\int_{\t^\th}^t d(X_r^{c},W_r)dr}e^{-c(t-\t^\th)}a_tdt\Big]\le \hE_{\tau^\th}\Big[\int_{\t^\th}^Te^{-\frac1\e\int_{\t^\th}^{t\wedge{\bar\t_\th^{c}}} d(X_r^c, W_r)dr}e^{-c(t-{\t^\th})}a_tdt\Big]\nonumber
 \\
 & \le& \hE_{\tau^\th  }\Big[\int_{\t^\th}^Te^{-\frac1\e \frac{\th}{2}[(t\wedge{\bar\t_\th^{c}})-\t^\th]}e^{-c(t-\t^\th)}a_tdt\Big]
 \le M \hE_{\tau^\th}\Big[\int_{\t^\th}^{\bar\t_\th^{c}}e^{-\frac\th{2\e}(t-\t^\th)}dt+\int_{\bar\t_\th^c}^Te^{-\frac\th{2\e}(\bar\t_\th^{c}-\t^\th)}dt\Big] \nonumber\\
 &\le& M\hE_{\t^\th}\Big[\int_{\t^\th}^Te^{-\frac\th{2\e }(t-\t^\th)}dt\Big]+M\hE_{\t^\th}[e^{-\frac\th{2\e}(\bar\t_\th^{c}-\t^\th)}]\dfnn A_\th(\e)+B_\th(\e),
 \eea
 where $A_\th(\cd)$ and $B_\th(\cd)$ are defined in an obvious way.
Clearly,  for fixed $\th=\th_0$,
\bea
\label{Ae0}
0\le A_\th(\e)\le\frac{2\e M}{\th}[1-e^{-\frac\th{2\e }T}]\to 0, \q \mbox{as $\e\to 0$, ~$\hP$-a.s.}
\eea
and the limit is uniform in  $(s,x,w)$ and $\pi\in\sU^{s,w}_{ad}[s,T]$.
 We shall argue that $B_\th(\e) \to 0$,  as $\e \to 0 $, in the same manner. Indeed, note that $X_{\t^\th}\le -\th$,  for $\d>0$ we have
 \bea
 \label{Ptth}
 \hP_{\t^\th}(|\bar\tau_\th^{c}-\tau^\th|<\delta)
& \le& \hP_{\t^\th}\Big\{\sup_{\tau^\th \le t \le\tau^\th+\delta}X_t^{\pi,\th,c}>-\frac{\th}{2}\Big\}
\le \hP_{\t^\th}\Big\{\sup_{\tau^\th \le t \le \tau^\th+\delta}[X_t^{\pi,\th,c}-X^\pi_{\t^\th}]>\frac{\th}{2}\Big\}\nonumber\\
&\le& \frac4{\th^2}\hE_{\t^\th}\Big\{\sup_{\tau^\th \le t \le \tau^\th+\delta}|X_t^{\pi,\th,c}-X^\pi_{\t^\th}|^2\Big\}\le C_\th\delta,
 \eea
 for some generic constant $C_\th>0$ depending only on $p$, $r$, $\si$, $T$,  $M$, and $\th$.
 Here we have applied Chebyshev inequality, as well as some standard SDE estimats.
Consequently, we derive from (\ref{Ptth}) that $ \sup_{\pi}\hP_{\t^\th}(|\bar\tau_\th^{c}-\tau^\th|<\delta)\le C_\th\delta$, $\hP$-a.s., and thus for fixed $\th$, and any $\eta>0$, we can find $\delta_0(\eta,\th)>0$, such that $\hP_{\t^\th}
(|\bar\tau_\th^{c}-\tau^\th|<\delta_0)<\frac{\eta}
{2T}$.
 Then,
 \bea
 B_\th(\e)&=&M
\Big\{\hE_{\t^\th}\Big[e^{-\frac{\th}{2\e}(\bar\tau_\th^{c}-\tau^\th)}:\bar\tau_\th^{c}-\tau^\th \geq \d_0\Big]+
 \hE_{\t^\th}\Big[e^{-\frac{\th}{2\e}(\bar\tau_\th^{c}-\tau^\th)}:\bar\tau_\th^{c}-\tau^\th < \delta_0\Big]\Big\} \\
 &\le&M\big\{ e^{-\frac{\th}{2\e }\delta_0}+\hP_{\t^\th}(\bar\tau_\th^{c}-\tau^\th< \delta_0)\big\}<Me^{-\frac{\th}{2\e}\delta_0}+\frac{\eta}2. \nonumber
 \eea
Therefore, for fixed $\th=\th_0$, one has $\overline\lim_{\e\to 0}B_\th(\e)\le \frac{\eta}2$, $\hP$-a.s. This, together with
(\ref{Jep}) and (\ref{Ae0}), then implies that
$\overline\lim_{\e \to 0} J^{\e}(\t^\th,X_{\t^\th}^{\pi},W_{\t^\th};\pi)\le \frac{\eta}2$, uniformly in $(s,x,w)\in K$ and $\pi\in\sU_{ad}[s,T]$, which in turn implies that, for $\th=\th_0$, $\overline\lim_{\e\to0}h^{\th}(\e)= \overline\lim_{\e\to0}\hE_{sxw}[V^\e(\tau^\th,X_{\tau^\th}^{\pi},W_{\tau^\th})]
\le \frac{\eta}2$, and the limit is uniformly in $(s,x,w)\in K$. Combining this with (\ref{Etht}) we derive from (\ref{VleVe1}) that
$$ V(s,x,w)\le \liminf_{\e\to0}V^\e(s,x,w)\le \limsup_{\e\to0} V^\e(s,x,w)\le V(s,x,w)+\eta.
$$
Since $\eta$ is arbitrary, we have $\lim_{\e\to 0}V^\e(s,x,w)=V(s,x,w)$, uniformly in $(s,x,w)\in K$.
Finally,  note that $V^\e$ is continuous in $x$, uniformly in $(s,x,w)\in K$, thanks to Lemma \ref{Vecont},
thus so is $V$. In particular,
$V$ is continuous in $x$ for $x\in[0,k]$, for all $k>0$,  proving the Theorem.
 \qed

\section{Continuity of the value function on $w$}
\setcounter{equation}{0}

We now turn our attention to the continuity of value function $V$ in the variable $w$. We should note that this is the most technical
part of the paper, as it involves the study of the delayed renewal process that has not been fully explored in the literature. We begin by the following proposition that extends the result of Proposition \ref{conts}. Recall the intensity of the interclaim times $T_i$'s:
$\l(t)=\frac{f(t)}{\bar F(t)}$, $t\ge 0$.
\begin{prop}
\label{conw1}
Assume that Assumption \ref{assump1} is in force. Then, for any  $h>0$  such that  $0\leq s<s+h<T$,
it holds that

(i) $V(s+h,x,w+h)-V(s,x,w)\leq\big[1-e^{-(ch+\int_{w}^{w+h}\l(u)du)}\big]V(s+h,x,w+h)$;

(ii) $V(s,x,w+h)-V(s,x,w) \leq Mh+\big[1-e^{-(ch+\int_{w}^{w+h}\l(u)du)}\big]V(s+h,x,w+h). $
\end{prop}

{\it Proof.}  (i)
For any $\pi=(\g, a)\in {\sU}^{s+h, w+h}_{ad}[s+h,T]$, we define, for $t\in[s,T]$,
$\tilde{\pi}^h_t=(\tilde\g_t, \tilde a_t)$ by
 \bea
\label{tildepi}
(\tilde{\gamma}_{t},\tilde{a}_{t})=(0,(p+r X^h_{t})\wedge M) )+[(\gamma_{t},a_{t})-(0,(p+r X^h_{t})\wedge M) )]{\bf 1}_{\{T_{1}^{s,w}> h\}}{\bf 1}_{[s+h,T]}(t).
\eea
	where $T_1^{s,w}$ is the first jump time of the delayed renewal process $N^{s,w}$, and
	$X^h:=X^{\tilde{\pi}^h, s,x,w}$. Since $T^{s,w}_1$ is a $\{\cF^s_t\}_{t\ge 0}=\{\cF_{s+t}\}_{t\ge0}$-stopping time, it is clear
	that $\tilde{\pi}^h\in \sU^{s,w}_{ad}[s,T]$. Let us  denote $\t^h:=\tau^{\tilde{\pi}^h}_{s,x,w}$ and
consider the following two cases:
	
	{\it Case 1.} $x\le\frac{M-p}{r}$.  In this case, for $s\le t <s+T^{s,w}_1$, we have $X^h_t\equiv x$ and $\tilde a_t\equiv p+rx\le M$. In particular, we note that by definition of $\tilde\pi^h$, given $T^{s,w}_1>h$ it must hold that $X^h_{s+h}=x$,
	$W^{s,w}_{s+h}=w+h$, and $T^{s+h,w+h}_1=T^{s,w}_1$, $\hP_{sxw}$-a.s. Thus
\bea
\label{Vest2}
V(s,x,w)&\ge &J(s,x,w;\tilde{\pi}^h)
\ge\hE_{sxw}\Big[\int_{s}^{\tau^{\tilde{\pi}^h}\wedge T}e^{-c (t-s)}\tilde{a}_{t}dt|T_1^{s,w}>h\Big]\hP_{sxw}\{T_1^{s,w}>h\}\nonumber\\
	&\ge&e^{-\int_{w}^{w+h}\l(u)du}\hE_{sxw}\Big[\int_{s+h}^{\tau^{\tilde{\pi}^h}\wedge T}e^{-c (t-s)}\tilde{a}_{t}dt
	\,\big|\,T_1^{s,w}>h\Big]\\
	&= &e^{-(ch+\int_{w}^{w+h}\l(u)du)}\hE_{(s+h)x(w+h)}\Big[\int_{s+h}^{\tau^{{\pi}}\wedge T}e^{-c (t-s-h)}{a}_{t}dt\Big]\nonumber\\
	&= &e^{-(ch+\int_{w}^{w+h}\l(u)du)}J(s+h,x,w+h;\pi). \nonumber
\eea
Since $\pi\in\sU_{ad}[s+h,T]$ is arbitrary, we obtain that $V(s,x,w)\geq e^{-(ch+\int_{w}^{w+h}\l(u)du)}V(s+h,x,w+h)$ which, with an  argument similar to the one led to (\ref{Vest1}), implies (a).
	
	{\it Case 2.} $x> \frac{M-p}{r}$,  In this case we have $\tilde a_s=M <p+r x=p+rX^h_s$, thus, by (\ref{Xsol}) $dX^{h}_s>0$.
	Namely, on the set $\{T^{s,w}_1>h\}$, $X^h$  will be continuous and increasing, so that $X^h_{s+h}=e^{rh}x+\frac{p-M}{r}(1-e^{-rh})=: x(h)$ (see (\ref{Xpi0})). Thus, noting that $W^{s,w}_{s+h}=w+h$
and $T^{s+h, w+h}_1=T^{s,w}_1$ on $\{T^{s,w}_1 >h\}$, a similar argument as (\ref{Vest2})  would lead to that
	$$V(s,x,w)\ge J(s,x,w;\tilde{\pi}^h)\ge e^{-(ch+\int_{w}^{w+h}\l(u)du)}V(s+h,x(h),w+h).
	$$
%
%
Now note that $x(h)>x$, it follows from Proposition \ref{vb}-(a) that $V(s+h,x(h),w+h))\ge V(s+h,x,w+h)$,
proving (a) again.
	
	Finally, (ii) follows from (i) and Proposition \ref{conts}-(b). This completes the proof.
\qed

The next result concerns the uniform continuity of $V$ on the variables $(s,w)$.
We have the following result.
\begin{prop}
\label{conw2}
Assume  that Assumption \ref{assump1} is in force. Then,
it holds that
$$\lim_{h\da 0}[V(s+h,x,w+h)-V(s,x,w)]=0,  \qq \mbox{uniformly in $(s,x,w)\in  D$.}
$$
\end{prop}

 {\it Proof.} From Proposition \ref{conw1}-(i) and the boundedness of $V$ we see that
\bea
\label{limit1}
\limsup_{h\da 0}[V(s+h,x,w+h)-V(s,x,w)] \le 0,  \qq \mbox{uniformly in $(s,x,w)\in D$.}
\eea
We need only prove the opposite inequality. We shall keep all the notations as in the previous proposition. For any $h\in(0, T-s)$,  and $\pi=(\gamma_{t}, a_{t})\in {\sU}_{ad}[s,T]$, we still consider  the strategy 
$\tilde{\pi}^h\in {\sU}^{s,w}_{ad}[s,T]$ defined by (\ref{tildepi}). (Note that $\tilde\pi^h$ depends on $\pi$ only for $t\in[s+h,T]$.) We again consider two cases, and denote
$\t_1:=T^{s,w}_1$ for simplicity.

 {\it Case 1. $x\le \frac{M-p}{r}$.}
In this case, we first write
	\begin{eqnarray}
	\label{82601}
	J(s,x,w;\tilde{\pi}^h)
	&=& \hE_{sxw}\Big[\int_{s}^{s+h}e^{-c (t-s)}\tilde{a}_{t}dt\Big|\t_1>h\Big]\hP(\t_1>h)\nonumber\\
	&&+ \hE_{sxw}\Big[\int_{s+h}^{\tau^{h}\wedge T}e^{-c (t-s)}\tilde{a}_{t}dt\Big|\t_1>h\Big]\hP(\t_1>h)\\
	&&+ \hE_{sxw}\Big[\int_{s}^{\tau^{h}\wedge T}e^{-c (t-s)}\tilde{a}_{t}dt\Big|\t_1\leq h\Big]\hP(\t_1\leq h):= I_{1}+I_{2}+I_{3},\nonumber
		\end{eqnarray}
	where $I_1, I_2$ and $I_3$ are defined as the three terms on the right hand side above, respectively. It is easy to see, by (\ref{tildepi}), that  on the set $\{\t_1>h\}$, $\tilde\g\equiv 0$, $X^{h}_t= x$, and $\tilde a_t= p+rx\le M$   for $t\in[s,s+h]$,  
	thus
	\begin{eqnarray}
	\label{82602}
	I_{1}
	&=&e^{-\int_{w}^{w+h}\lambda(u)du}\hE_{sxw}\Big[\int_{s}^{s+h}e^{-c (t-s)}(p+rX^{h}_{t})dt|\t_1>h\Big]
	\leq (p+rx)h;\\
	I_{2}
	&\le& e^{-ch-\int_{w}^{w+h}\lambda(u)du}V(s+h,x,w+h)\le V(s+h,x,w+h). \nonumber
	\label{82603}
	\end{eqnarray}
To estimate $I_3$, we first note that on the set $\{\t_1\le h\}$, by (\ref{tildepi}), $\tilde\g_t\equiv 0$, for all $t\in[s,T]$.
Thus $X^{h}_t= x$ and $\tilde a_t= p+rx$ for $t\in[s, s+\t_1)$. We also note that $\t^{h}\ge s+\t_1$ and  $\{\t^{h}>s+\t_1\}=\{U_1\le x\}$. Bearing these in mind we now write
	\bea
	\label{I3}
	I_{3}= 	\hE_{sxw}\Big[\Big(\int_s^{s+\t_1}+\int_{s+\t_1}^{\tau^{h}\wedge T}\Big)e^{-c (t-s)}\tilde{a}_{t}dt:\t_1\leq h\Big]:=I^1_3+I^2_3,
	\eea
where $I^1_3$ and $I^2_3$ are defined in an obvious way. For simplicity let us denote the density function of $T^{s,w}_1$ by
$p_{\t_1}(z)=\lambda(w+z)e^{-\int_{w}^{w+z}\lambda(v)dv}$, $z\ge0$. Clearly, given $\t_1\le h$ we have
\bea
\label{I31}
I^1_3	&=&\int_{0}^{h}\hE_{sxw}\Big[\int_{s}^{s+\t_1}e^{-c (t-s)}(p+rX^{h}_{t})dt|\t_1=z\Big]p_{\t_1}(z)dz \nonumber\\
&=&\int_{0}^{h}\Big[\int_{s}^{s+z}e^{-c (t-s)}(p+rx)dt\Big]p_{\t_1}(z)dz\\
&\le& \int_s^{s+h}e^{-c (t-s)}(p+rx)dt (1-e^{-\int_{w}^{w+h}\lambda(v)dv})\le (1-e^{-\int_{w}^{w+h}\lambda(v)dv})(p+rx)h. \nonumber
\eea
Further, we note that $(X^{h}_{s+\t_1},W^{s,w}_{s+\t_1})=(x-U_1,0)$,
 $\hP$-a.s., thus
\bea
\label{I32}
I^2_3&=& \int_{0}^{h}\hE_{sxw}\Big[\int_{s+z}^{\tau^{h}\wedge T}e^{-c (t-s)}(p+rX^{h}_{t})dt{\bf 1}_{\{\t^{h}>s+z\}}|\t_1=z\Big]p_{\t_1}(z)dz \nonumber\\
	&=&
\int_{0}^{h}\int_{0}^{x}\hE_{sxw}[\int_{s+z}^{\tau^{h}\wedge T}e^{-c (t-s)} (p+rX^{h}_{t})dt|\t_1=z, U_{1}=u]p_{\t_1}(z)dG(u) dz\\
	&\leq& \int_{0}^{h}\int_{0}^{x}e^{-cz}V(s+z,x-u,0) 	p_{\t_1}(z)dG(u) dz\le \frac{M}c (1-e^{-\int_{w}^{w+h}\lambda(v)dv}).\nonumber
	\eea
Here the last inequality is due to Proposition \ref{vb}-(ii). Now, combining (\ref{I31}) and (\ref{I32}) we have
	\bea
 \label{82604}
I_3\le
(1-e^{-\int_{w}^{w+h}\lambda(v)dv})((p+rx)h+\frac{M}c),
	\eea
and consequently we obtain from (\ref{82601})-(\ref{82604}) that, for $x<\frac{M-p}{r}$,
	\begin{eqnarray}
	\label{Jest1}
	J(s,x,w;\tilde{\pi}^h)\leq (p+rx)h+V(s+h,x,w+h)+(1-e^{-\int_{w}^{w+h}\lambda(v)dv})((p+rx)h+M/c).
	\end{eqnarray}

{\it Case 2. $x\ge\frac{M-p}{r}$.} In this case, using the strategy $\tilde\pi^h$ as in (\ref{tildepi}) with a similar argument
as in Case 1 we can derive that
\begin{eqnarray}
\label{Jest2}
J(s,x,w;\tilde{\pi}^h)&\leq& Mh+V(s+h,e^{r h}(x+\frac{p-M}{r}(1-e^{-rh})), w+h)\nonumber\\
&&+(1-e^{-\int_{w}^{w+h}\lambda(v)dv})(M(h+\frac{1}c)).
\end{eqnarray}

To complete the proof we are to replace the left hand side of (\ref{Jest1}) and (\ref{Jest2}) by $J(s,x,w, \pi)$, which would lead to
the desired inequality, as $\pi\in {\cal U}_{ad}[s,T]$ is arbitrary. To this end we shall
argue along a similar line as those in the previous section.

Recall the penalty function $\beta^{\pi,s}(t, \e):=
\beta^{\pi, s,x,w}(t,\e)$ defined by (\ref{beta}), and define
$$J^{\e}(s,w,x;\pi)=\hE_{swx}\Big[\int_s^{T}\beta^{\pi,s}(t,\e) e^{-c(t-s)}a_tdt\Big].
$$
We first write
\bea
&&\left|J^{\e}(s,x,w;\pi)-J^{\e}(s,x,w;\tilde\pi^h)\right|\nonumber\\
&\le& \hE_{sxw}\left|\int_s^{s+h}e^{-c(t-s)}[\beta^{\pi,s}(t, \e)a_t-\beta^{\tilde\pi^h,s}(t, \e)\tilde{a_t}]dt\right|\\
&&+\hE_{sxw}\left|\int_{s+h}^Te^{-c(t-s)}[\beta^{\pi,s}(t,\e)a_t-\beta^{\tilde\pi^h,s}(t, \e)\tilde a_t]dt\right|: = I_1+I_2\notag
\eea
It is easy to see that $I_1<2Mh$, thanks to Assumption \ref{assump1}.
We shall estimate $I_2$. Note that
\bea
\label{I2}
I_2&=&\hE_{sxw}\Big\{\Big|\int_{s+h}^Te^{-c(t-s)}(\beta^{\pi,s}(t,\e)-\beta^{\tilde\pi^h,s}(t,\e))a_tdt\Big|\Big|\t_1>h\Big\}\hP(\t_1>h)\\
&&+\hE_{sxw}\Big\{\Big|\int_{s+h}^Te^{-c(t-s)}[\beta^{\pi,s}(t,\e)a_t-\beta^{\tilde\pi^h,s}(t,\e))\tilde a_t]dt\Big|\Big|\t_1\le h\Big\}\hP(\t_1\le h):= I^1_2+I^2_2.\notag
\eea
Since $X^\pi_t \ge 0$, $X^h_t\ge0$ for $t\le s+h$ on the set $\t_1>h$ (i.e., ruin occurs only at arrival of a claim), we have
$d(X_t^{\pi},W_t)=d(X_t^{h},W_t)=0$ for $t\in[s,s+h]$, i.e., $\beta^{\pi,s}(t,\e)=\beta^{\pi,s+h}(t,\e)$,
$\beta^{\tilde\pi^h,s}(t,\e)=\beta^{\tilde\pi^h,s+h}(t,\e)$, for $t\in[s+h,T]$. Thus, by the similar arguments as in Lemma \ref{Vecont} one
shows that
\bea
\label{I21}
I_2^1&=&\hE_{sxw}\Big\{\Big|\int_{s+h}^T(\beta^{\pi,s+h}(\e,t)-\beta^{\tilde\pi^h,s+h}(\e,t))e^{-c(t-s)}a_tdt\Big|\Big|\t_1>h\Big\}\hP(\t_1>h)\nonumber\\
&\le& C\hE_{sxw}|X_{s+h}^{\pi}-X_{s+h}^{h}|,
\eea
where $C>0$ is a generic constant depending only on $\e$ and $T$. Furthermore, since $\hP(\t_1\le h)=(1-e^{-\int_w^{w+h}\l(v)dv})=O(h)$, we have $I_2^2=O(h)$. It then follows from (\ref {I2}) and (\ref{I21}) that
$I_2 \le C\hE_{sxw}|X_{s+h}^{\pi}-X_{s+h}^{h}|+O(h)
$.
The standard result of SDE then leads to $\lim_{h\to 0}I_2=0$, whence
$\lim_{h\to0}|J^{\e}(s,x,w;\pi)-J^{\e}(s,x,w;\tilde\pi^h)|= 0$, and the convergence is obviously uniform for $(s,x,w)\in D$ and
$\pi\in\sU^{s,w}_{ad}[s,T]$.

To complete the proof we note that, with exactly the same argument as that in Theorem \ref{continuityx} one shows that,
for any $\eta>0$, there exists $\e_0>0$, such that
$$|J^{\e_0}(s,x,w;\pi)-J(s,x,w;\pi)|+|J^{\e_0}(s,x,w;\tilde\pi^h)-J(s,x,w;\tilde\pi^h)|<\eta, \q\forall (s,x,w)\in D.
$$
Then, for the fixed $\e_0$, we choose $h_0>0$, independent of $\pi\in\sU^{s,w}_{ad}[s,T]$ such that
$$ |J^{\e_0}(s,x,w;\pi)-J^{\e_0}(s,x,w;\tilde\pi^h)|<\eta, \q \forall (s,x,w)\in D, \q \forall 0<h<h_0.
$$
Thus, if $x<\frac{M-p}{r}$, for all $0<h<h_0$, we derive from (\ref{Jest1}) that
\beaa
&& J(s,x,w;\pi)-V(s+h,x,w+h)\notag\\
& \le& |J(s,x,w;\pi)-J^{\e_0}(s,x,w;\pi)|+|J^{\e_0}(s,x,w;\pi)-J^{\e_0}(s,x,w;\tilde\pi^h)|\\
&&+|J^{\e_0}(s,x,w;\tilde\pi^h)-J(s,x,w;\tilde{\pi}^h)|+J(s,x,w;\tilde{\pi}^h)-V(s+h,x,w+h)\notag\\
& \le& 2\eta+(p+rx)h+(1-e^{-\int_{w}^{w+h}\l(v)dv})((p+rx)h+M/c) \le 2\eta + g_1(h).
\nonumber
\eeaa
 where $g_1(h):=Mh+(1-e^{-\int_{w}^{w+h}\l(v)dv})(Mh+M/c)$.
Since $\pi\in\sU_{ad}[s, T]$ is arbitrary, we have
\bea
\label{xleMpr} V(s, x,w)-V(s+h, x, w+h)\le  2\eta+g_1(h).
\eea
First sending $h\to0$ and then $\eta\to0$ we obtain the desired opposite inequality of (\ref{limit1}).

The case for $x\ge\frac{M-p}{r}$ can be argued similarly. We apply (\ref{Jest2}) to get the analogue of (\ref{xleMpr}):
\bea
\label{xgeMpr}
&& V(s,x,w)-V(s+h,x,w+h)\\
& \le& 2\eta +g_1(h)+V(s+h,e^{r h}(x+\frac{p-M}{r}(1-e^{-rh})), w+h)-V(s+h,x,w+h).\notag
\eea
For fixed $x\geq \frac{M-p}{r}$, by first sending $h\to0$ and then
$\eta\to0$, we have
\bea
\liminf_{h\da 0}[V(s+h,x,w+h)-V(s,x,w)] \geq 0.
\eea
thanks to the uniformly continuity  $V(s,x,w)$ in $x$ (uniformly in $(s,w)$). This, together with (\ref{limit1}), yields that, for given $x\geq 0$,
\bea
\label{32401}
\lim_{h\da 0}[V(s+h,x,w+h)-V(s,x,w)]=0,  \qq \mbox{uniformly in $(s,w)$.}
\eea
Then, combining (\ref{32401}) and Proposition \ref{conw1}, one shows that $V(s,x,w)$ is continuous in $(s,w)$ for fixed $x$.
It remains to argue that (\ref{32401}) holds uniformly in $(s,x,w)\in D$.

To this end, we note that, by Proposition \ref{vb} and Theorem  \ref{continuityx},  $V(s,x,w)$ is increasing in $x$, continuous in $(s,w)$, and with
a continuous limit function $\frac{M}{c}(1-e^{-(T-s)})$ (in $(s,w)$). Thus $V(s,x,w)$ converges uniformly to  $\frac{M}{c}(1-e^{-(T-s)})$ as $x\rightarrow \infty$, uniformly in $(s,w)$, thanks to Dini's Theorem. That is,
for $\eta>0$, there exists
$N=N(\eta)>\frac{M-p}{r}$,  such that
$$V(s+h,e^{r h}(x+\frac{p-M}{r}(1-e^{-rh})), w+h)-V(s+h,x,w+h)< \eta, \qq x>N.$$
On the other hand, for $\frac{M-p}{r}\leq x \leq N$, by Theorem \ref{continuityx},  there exists $\delta(\eta)=\d(N(\eta))>0$,
such that for $h<\delta(N)$, it holds that
$$V(s+h,e^{r h}(x+\frac{p-M}{r}(1-e^{-rh})), w+h)-V(s+h,x,w+h)< \eta.$$
Thus, we see from (\ref{xgeMpr}) that  for all $(s,x,w)\in D$, and $x\ge \frac{M-p}r$,
$$ V(s,x,w)-V(s+h,x,w+h) \le 4\eta, \qq \mbox{whenever  $h<\delta$.} $$
Combining this with the case $x<\frac{M-p}r$ argued previously, we see that
$$\liminf_{h\da 0}[V(s+h,x,w+h)-V(s,x,w)]\ge 0,  \qq \mbox{uniformly in $(s,x,w)\in D$,}
$$
proving the opposite inequality of (\ref{limit1}), whence the proposition.
\qed

Combining Theorems \ref{conts} and \ref{conw1}, we have proved the following theorem.
\begin{thm}
\label{coninw}
Assume that Assumption \ref{assump1} is in force. Then, the value function $V(s,x,w)$ is uniformly continuous in $w$, uniformly on $(s,x,w)\in D$.
\qed
\end{thm}


\section{Dynamic Programming Principle}
\setcounter{equation}{0}

In this section we shall substantiate the Bellman Dynamic Programming Principle (DPP) for our optimization problem. We begin
with a simple but important lemma.
	\begin{lem}
\label{DPPlemma}
For any $\varepsilon>0$, there exists $\delta>0$,  independent of $(s,x,w)\in D$, such that for any $\pi\in \sU_{ad}^{s,w}[s,T]$
 and $h:=(h_{1}, h_{2})$ with $0\le h_{1}, h_{2}<\delta$, we can find $\hat{\pi}^h \in\sU^{s,w-h_2}_{ad}[s,T]$ such that
\bea
\label{Jep12}
J(s,x,w,\pi)-J(s,x-h_{1},w-h_{2},\hat{\pi}^{h})\leq \varepsilon, \q \forall (s,x,w)\in D.
\eea
Moreover, the construction of $\hat\pi^h$ is independent of $(s,x,w)$.
\end{lem}

{\it Proof.} Let $\pi=(\g, a)\in \sU^{s,w}_{ad}[s,T]$. For any $h=(h_1, h_2)\in [0,\infty)^2$, we consider the following two modified strategies in the form of (\ref{tildepi}): denoting $\th(x):=(p+rx)\wedge M$,
\bea
\label{tildehatpi}
\left\{\ba{lllll}
\tilde\pi^{h}_t:=(\tilde\gamma_{t}^{h},\tilde a_{t}^{h})=(0,\th(\tilde X^{h}_{t}))+[(\gamma_{t},a_{t})
-(0,\th(\tilde X^{h}_{t})]{\bf 1}_{\{\tilde \t_{1}^{h}>h_{2}\}}{\bf 1}_{[s, T]}(t), \qq t\in [s-h_{2}, T]; \bs\\
\hat \pi^{h}_t:=(\hat\gamma_{t}^{h},\hat a_{t}^{h})=(0,\th(\hat X^{h}_{t}))+[(\gamma_{t-h_2},a_{t-h_2})-
(0,\th(\hat X^{h}_{t}))]{\bf 1}_{\{\hat \t_{1}^{h}>h_{2}\}}{\bf 1}_{[s+h_{2},T]}(t),  ~t\in [s, T].
\ea\right.
\eea
where, for notational simplicity, we denote $\tilde \t^h_1:=T^{s-h_2, w-h_2}_1$; $\hat \t^h_1:=T_{1}^{s,w-h_2}$; $\tilde X^{h}:=X^{ \tilde \pi^{h},s-h_{2},x,w-h_{2}}$; and $ \hat X^{h}:=X^{\hat\pi^{h},s,x,w-h_{2}}$. Clearly, $\tilde\pi^{h}\in \sU^{s-h_2,w-h_2}_{ad}[s-h_{2},T]$ and  $\hat\pi^{h}\in \sU^{s,w-h_2}_{ad}[s,T]$, and it holds that
\beaa
&&J(s,x,w;\pi)-J(s,x-h_{1},w-h_{2};\hat\pi^{h})\\
&\leq& [J(s,x,w,\pi)-J(s-h_{2},x,w-h_{2};\tilde\pi^{h})]
+[J(s-h_{2},x,w-h_{2};\tilde\pi^{h})-J(s,x,w-h_{2},\hat \pi^{h})]\\
&&+[J(s,x,w-h_{2},\hat\pi^{h}) -J(s,x-h_{1},w-h_{2},\hat\pi^{h})]:= J_1+J_2+J_3.
\eeaa
We shall estimate $J^i$'s separately. First, by (\ref{Vest2}), we have
\bea
\label{J1}
J_1&=&J(s,x,w,\pi)-J(s-h_{2},x,w-h_{2},\tilde\pi^{h})\leq  [1-e^{-(c h_{2}+\int_{w-h_2}^{w}\l(u)du)}\big]J(s,x,w,\pi)\nonumber\\
&\leq & \frac{M}{c}[1-e^{-(c h_{2}+\int_{w-h_{2}}^{w}\l(u)du)}\big].
\eea
Next, we observe from definition (\ref{tildehatpi}) that the law of $\tilde X^{h}$ on $[s-h_2, T-h_2]$ and that of $\hat X^{h}$
on $[s, T]$ are identical. We have
\bea
\label{J20}
J_2&=&J(s-h_{2},x,w-h_{2},\tilde\pi^{h})-J(s,x,w-h_{2},\hat\pi^{h})\nonumber\\
&=&\hE_{(s-h_2)x(w-h_2)}\Big[\int_{s-h_{2}}^{\t^{\tilde\pi^h}\wedge T}e^{-c(t-s+h_2)}\tilde a_{t}^{h}dt\Big]-\hE_{sx(w-h_2)}\Big[\int_{s}^{\tau^{\hat\pi^h}\wedge T}e^{-c(t-s)}\hat a_{t}^{h}dt\Big]\\
&=&e^{-ch_{2}}\hE_{(s-h_2)x(w-h_2)}\Big[\int_{s-h_{2}}^{\tau^{\tilde\pi^h}\wedge (T-h_{2})}e^{-c(t-s)}\tilde a_{t}^{h}dt\Big]-\hE_{sx(w-h_2)}\Big[\int_{s}^{\tau^{\hat\pi^{h}}\wedge T}e^{-c(t-s)}\hat a_{t}^{h}dt\Big] \nonumber\\
&&+ \hE_{(s-h_2)x(w-h_2)}\Big[\int_{\tau^{\tilde\pi^h}\wedge (T-h_{2})}^{\t^{\tilde\pi^h}\wedge T}e^{-c(t-s+h_{2})}\tilde a_{t}^{h}dt\Big]:= e^{-ch_2}J^1_2-J^2_2+J^3_2,  \nonumber
\eea
where $J^i_2$, $i=1,2,3$ are the three expectations on the right side, respectively. Note that by definition of the $\hat \pi^h$ and $\tilde \pi^h$, it is easy to check that $J^1_2=J^2_2$. Thus (\ref{J20}) becomes
\bea
\label{J2}
J_2\le J^3_2=\hE_{(s-h_2)x(w-h_2)}\Big[\int_{\tau^{\tilde\pi^h}\wedge (T-h_{2})}^{\t^{\tilde\pi^h}\wedge T}e^{-c(t-s+h_{2})}\tilde a_{t}^{h}dt\Big]\leq Mh_2. 
\eea
Finally, from the proofs of  Theorem \ref{continuityx} and Lemma \ref{Vecont}, we see that the mapping  $x\mapsto J(s,x,w, \pi)$ is continuous in $x$, uniformly for $(s,x,w)\in D$ and $\pi\in \sU_{ad}[s, T]$. Therefore, for any $\e>0$, we can find $\d>0$, depending only on $\e$, such that, for $0<h_1<\d$, it holds that
\beaa
J_3= J(s,x,w-h_{2},\hat\pi^{h}) -J(s,x-h_{1},w-h_{2},\hat\pi^{h})<\e/3, \qq \mbox{$\forall h_2\in(0,w)$.}
\eeaa
We can then assume that $\d$ is small enough, so that for $h_2<\d$, it holds that  $J_1<\e/3$, $J_2<\e/3$, uniformly in
$(s,x,w)\in D$ and $\pi\in \sU_{ad}[s,T]$, thanks to (\ref{J1}) and (\ref{J2}). Consequently, we have
\beaa
J(s,x,w,\pi)-J(s,x-h_{1},w-h_{2},\hat\pi^{h})
%
\le J_1+J_2+J_2<\varepsilon,
\eeaa
proving (\ref{Jep12}), whence the lemma.
\qed

We are now ready to prove the first main result of this paper: the Bellman Principle of Optimality or Dynamic Programming
Principle (DPP). Recall that for a given $\pi\in\sU_{ad}[s,T]$ and $(s,x,w)\in D$, we denote $R^\pi_t=R^{\pi, s,x,w}_t=
(t, X^{\pi, s,x,w}_t, W^{s,w}_t)$, $t\in[ s,T]$.
\begin{thm}
\label{DPPthm}
Assume that Assumption \ref{assump1} is in force. Then, for any $(s, x,w)\in D$ and any stopping time $\t\in [ s, T]$,
it holds that
\bea
\label{DPP}
V(s,x,w) = \mathop {\sup }\limits_{\pi  \in {{\sU} _{ad}[s,T]}} {\hE_{sxw}}\Big[\int_s^{\t  \wedge {\t ^\pi }} {{e^{ - c(t - s)}}a_tdt + {e^{ - c(\t  \wedge {\t ^\pi } - s)}}V({R^{\pi}_{\t  \wedge {\t ^\pi }}})} \Big].
\eea
\end{thm}

{\it Proof.} The idea of the proof is more or less standard. We shall first argue that (\ref{DPP}) holds for deterministic $\t=s+h$,
for $h\in(0,T-s)$. That is, denoting
\beaa
\label{DPPf}
v(s,x,w;s+h):=\sup_{\pi  \in {{\sU} _{ad}[s,T]}} \hE_{sxw}}\Big[\int_s^{(s+h) \wedge {\t ^\pi }} {{e^{ - c(t - s)}}a_tdt + {e^{ - c((s+h) \wedge {\t ^\pi } - s)}}V(R_{(s+h) \wedge {\t ^\pi }}^{\pi} )\Big],
\eeaa
we are to show that $V(s,x,w)=v(s,x,w; s+h)$. To this end, let
$\pi  = (\g,a)\in \sU_{ad}[s, T]$, and write
\bea
\label{vs+h}
J(s,x,w;\pi) =
\hE_{sxw}\Big[\int_s^{(s+h)\wedge \t^\pi} {{e^{ - c(t - s)}}a_tdt\Big] + \hE_{sxw}\Big[\int_{s+h}^{{\t ^\pi }} {{e^{ - c(t - s)}}a_tdt} } :\t^\pi>s+h\Big].
\eea
%
%
%
%
Now applying Lemma \ref{PDPP} we see that the second term on the right hand side of (\ref{vs+h}) becomes
\beaa
 &&\hE_{sxw}\Big[\int_{s+h}^{\t ^\pi }e^{ - c(t - s)}a_tdt:\t^\pi>s+h\Big]\nonumber\\
 &=&e^{-ch}\hE_{sxw}\Big[\hE\Big[\int_{s+h}^{\t ^\pi}e^{ - c(t - (s+h))}a_tdt\Big| {\cal F}_{s+h}^{s}\Big]:\t^\pi>s+h\Big]\\
&=&e^{-ch}\hE_{sxw}\Big[J^{\pi}(s+h, X^{\pi}_{s+h}, W^{\pi;}_{s+h}) :\t^\pi>s+h\Big]\le e^{ - ch}\hE_{sxw}
\Big[V(R_{s+h}^{\pi}):\t^\pi>s+h\Big] \nonumber\\
 &\le &\hE_{sxw}\Big[e^{ - c((s+h) \wedge {\tau ^\pi } - s)} V(R_{(s+h) \wedge {\tau ^\pi }}^{\pi})\Big].
\eeaa
Plugging this into (\ref{vs+h}) and taking supremum on both sides above we obtain that  $V(s,x,w) \le v(s,x,w;s+h)$.

The proof of the reversed inequality is slightly more involved, as usual.
To begin with, we recall Lemma \ref{DPPlemma}. For any $\e>0$, let $\d>0$ be the constant in Lemma \ref{DPPlemma}.
Next, let $0=x_0<x_1<\cds $ and $0=w_0<w_1<\cds <w_n=T$ be a partition of $[0,\infty)\times [0,T]$, so that $x_{i+1}-x_i<\d$
$w_{j+1}-w_j<\d$. Denote $D_{ij}:=[x_{i-1},x_i)\times[w_{j-1},w_j)$,  $i,j\in\hN$. For
$0\le s<s+h<T$, $i\in\hN$, and $0\le j\le n$ we choose $\pi^{ij}\in {{\sU}^{s+h, w_j}_{ad}[s+h,T]}$ such that
\beaa
J(s+h,x_i,w_j;\pi^{ij})>V(s+h,x_i,w_j)-\varepsilon.
\eeaa
Now applying Lemma $\ref{DPPlemma}$, for each $(x,w)\in D_{ij}$ and $\pi^{ij}\in \sU_{ad}^{s+h,w_j}[s+h, T]$, we can define strategy $\hat\pi^{ij}=\hat\pi^{ij}(x,w)\in \sU_{ad}^{s+h,w}[s+h,T]$, such that
\bea
\label{3e}
J(s+h,x,w;\hat\pi^{ij})&\ge&   J(s+h,x_i,w_j;\pi^{ij})-\varepsilon\notag\\
&\ge& V(s+h,x_i,w_j)-2\varepsilon \ge V(s+h,x,w)-3\varepsilon.
\eea
In the above the last inequality is due to the uniform continuity of $V$  on the variables $(x,w)$.

Now for any $\pi\in\sU_{ad}^{s,w}[s,T]$, we define a new strategy $\pi^*$ as follows:
\beaa
\label{pistar}
\pi^*_t=\pi_t{\bf 1}_{[s, s+h)}(t)+\sum_{i=0}^\infty\sum_{j=0}^{n-1} \hat\pi^{ij}_t(X^\pi_{s+h},W_{s+h}){\bf 1}_{D_{ij}}(X^\pi_{s+h}, W_{s+h}){\bf 1}_{[s+h, T]}(t).
\eeaa
%
%
%
Then one can check that $\pi^*\in \sU_{ad}^{s,w}[s, T]$, and  $\{\t^{\pi^*}\le s+h\}=\{\t^\pi\le s+h\}$. Furthermore,  
when $\t^{\pi^*}>s+h$ we have
\bea
J(s+h,X^\pi_{s+h},W_{s+h};\pi^*)\ge V(s+h,X^\pi_{s+h},W_{s+h})-3\e, \q \mbox{$\hP$-a.s. on $\{\t^{\pi^*}>s+h\}$},
\eea
thanks to (\ref{3e}). Consequently, similar to (\ref{vs+h}) we have
\bea
\label{Vgev}
V(s,x,w)&\ge& J(s,x,w;\pi^*)\\
&=&\hE_{sxw}\Big[ \int_s^{(s+h)\wedge\t^\pi} {{e^{ - c(t - s)}}a_tdt}  +{\bf 1}_{\{\t^{\pi}>s+h\}}e^{-ch}\int_{s+h}^{\tau ^{\pi ^*\wedge T}} e^{ - c(t - (s+h))}a_t^*dt\Big]\nonumber\\
&=&\hE_{sxw}\Big[\int_s^{(s+h)\wedge\t^\pi} {{e^{ - c(t - s)}}a_td{t}}  + {\bf 1}_{\{\t^{\pi}>s+h\}}e^{ - ch}J(s+h,{X^\pi_{s+h}},{W_{s+h}};\pi ^*)\Big] \notag\\
&\ge&\hE_{sxw}\Big[\int_s^{(s+h)\wedge\t^\pi} {{e^{ - c(t - s)}}a_td{t}}  + {e^{ - c({(s+h)\wedge\t^\pi} - s)}}V(R^\pi_{(s+h)\wedge\t^\pi})\Big]  - 3\varepsilon.\nonumber
\eea
Here in the last inequality we used the fact that ${\bf 1}_{\{\t^\pi\le s+h\}}V(R^\pi_{(s+h)\wedge\t^\pi})={\bf 1}_{\{\t^\pi\le s+h\}}V(R^\pi_{\t^\pi})=0$.
Since $\pi$ is arbitrary, (\ref{Vgev}) implies $V(s,x,w)\ge v(s,x,w;s+h)-3\e$. Since $\e>0$ is arbitrary, we obtain
that $V(s,x,w)\ge v(s,x,w;s+h)$, proving (\ref{DPP}) for $\t=s+h$.

We now consider the general case when $s<\t<T$ is a stopping time. Let $s=t_0<t_1<\cds<t_n=T$ be a partition of $[s,T]$.
We assume that $t_k:= s+\frac{k}{n}(T-s)$, $k=0,1,\cds, n$. Define
$\t_n:=\sum_{k=0}^{n-1} t_{k}{\bf 1}_{[t_{k}, t_{k+1})}(\t)$.
Clearly,  $\t_n$ takes only a finite number of values and $\t_{n}\to\t$, $\hP$-a.s. It is easy to check, using the same argument
above when $\t$ is deterministic to each subinterval $[s,  T]$, that $V(s,x,w)\le v(s,x,w; \t_n)$. We shall prove by induction
(on $n$) that
\bea
\label{Vgevn}
V(s,x,w)\ge v(s,x,w;\t_n), \qq \forall n\ge 1.
\eea
 Indeed,
for $n=1$, we have $\t_1\equiv s$, so there is nothing to prove. Now suppose that (\ref{Vgevn}) holds for $\t_{n-1}$, and
$n\ge 2$. We shall argue that  (\ref{Vgevn}) holds for $\t_n$ as well.
For any $\pi\in \sU_{ad}^{s,w}[s, T]$ we have
\bea
\label{vtn}
&&\hE_{sxw}\Big\{\int_s^{\t_n \wedge {\t ^\pi }} {{e^{ - c(t - s)}}a_td{t} + {e^{ - c(\t_n \wedge {\t ^\pi } - s)}}V(R_{\t_n \wedge {\t ^\pi }}^{\pi} )}\Big\}\notag\\
&=&\hE_{sxw}\Big\{{\bf 1}_{\{\t^{\pi}\le t_1\}}\int_s^{ {\t ^\pi }} {e^{ - c(t - s)}}a_td{t}\Big\}\\
&&+\hE_{sxw}\Big\{\Big[\int_s^{\t_n \wedge {\t ^\pi }} {{e^{ - c(t - s)}}a_td{t} + {e^{ - c(\t_n \wedge {\t ^\pi } - s)}}V(R_{\t_n \wedge {\t ^\pi }}^{\pi} )}\Big]{\bf 1}_{\{\t_n>t_1\}}{\bf 1}_{\{\t^{\pi}>t_1\}} \nonumber\\
&&+\Big[\int_s^{t_1 } {{e^{ - c(t - s)}}a_td{t} +e^{ - c(t_1 - s)}V(R_{t_1 }^{\pi} )}\Big]{\bf 1}_{\{\t_n=t_1\}}{\bf 1}_{\{\t^{\pi}>t_1\}}\Big\}.\notag
\eea

Note that on the set $\{\t_n>t_1\}$, $\t_n$  takes only  $n-1$ values, by inductional hypothesis,
we have
\beaa
&&\hE_{sxw}\Big\{\Big[\int_{t_1}^{\t_n \wedge {\t ^\pi }} {{e^{ - c(t - s)}}a_td{t} + {e^{ - c(\t_n \wedge {\t ^\pi } - s)}}V(R_{\t_n \wedge {\t ^\pi }}^{\pi} )}\Big]{\bf 1}_{\{\t_n>t_1\}}{\bf 1}_{\{\t^{\pi}>t_1\}}\Big\}\\
&\le &\hE_{sxw}\Big\{e^{-c(t_1-s)}v(t_1, X^\pi_{t_1}, W_{t_1};\t_n)
{\bf 1}_{\{\t_n>t_1\}}{\bf 1}_{\{\t^{\pi}>t_1\}}\Big\}\le\hE_{sxw}\Big\{e^{-c(t_1-s)}V(R^\pi_{t_1}){\bf 1}_{\{\t_n>t_1\}}{\bf 1}_{\{\t^{\pi}>t_1\}}\Big\}.
\eeaa
Plugging this into (\ref{vtn}) we obtain
\bea
&&\hE_{sxw}\Big\{\int_s^{\t_n \wedge {\t ^\pi }} {{e^{ - c(t - s)}}a_td{t} + {e^{ - c(\t_n \wedge {\t ^\pi } - s)}}V(R_{\t_n \wedge {\t ^\pi }}^{\pi} )}\Big\} \nonumber\\
&\le&\hE_{sxw}\Big\{{\bf 1}_{\{\t^{\pi}\le t_1\}}\int_s^{ {\t ^\pi }} {e^{ - c(t - s)}}a_td{t}\Big\}+ \hE_{sxw}\Big\{\Big[\int_s^{t_1 }e^{ - c(t - s)}a_td{t} + e^{ - c(t_1 - s)}V(R_{t_1 }^{\pi} )\Big]   {\bf 1}_{\{\t_n>t_1\}}{\bf 1}_{\{\t^{\pi}>t_1\}}
\notag\\
&&+ \Big[\int_s^{t_1 }e^{ - c(t - s)}a_td{t} + e^{ - c(t_1 - s)}V(R_{t_1 }^{\pi} )\Big]{\bf 1}_{\{\t_n=t_1\}}{\bf 1}_{\{\t^{\pi}>t_1\}}\Big\}\notag\\
&=& \hE_{sxw}\Big\{{\bf 1}_{\{\t^{\pi}\le t_1\}}\int_s^{ {\t ^\pi }} {e^{ - c(t - s)}}a_td{t}\Big\}+ \hE_{sxw}\Big\{{\bf 1}_{\{\t^{\pi}>t_1\}}
\Big[e^{-c(t_1-s)}V(R_{t_1}^{\pi})+ \int_s^{t_1 } {e^{ - c(t - s)}}a_t{dt}\Big]\Big\}\notag\\
&=& \hE_{sxw}\Big\{\int_s^{ t_1\wedge {\t ^\pi }} {e^{ - c(t - s)}}a_td{t}+ e^{-c(t_1\wedge \t^\pi-s)}V(R_{t_1\wedge\t^\pi}^{\pi})\Big\}\le V(s,x,w).  \nonumber
\eea
In the above we again used the fact $V(R^\pi_{\t^\pi})=0$, and the last inequality is due to (\ref{DPP}) for fixed time $t_1=s+h$.
Consequently we obtain $v(s,x,w;\t_n)\le V(s,x,w)$, whence $v(s,x,w;\t_n)=V(s,x,w)$. A simple application of Dominated
Convergence Theorem, together with the uniform continuity of the value function, will then leads to the general form of (\ref{DPP}).
The proof is now complete.
\qed

\section{The Hamilton-Jacobi-Bellman equation.}
\setcounter{equation}{0}
	
	We are now ready to investigate the main subject of the paper: the Hamilton-Jacobi-Bellman (HJB) equation associated to
our optimization problem (\ref{V1}). We note that such a PDE characterization of the value function is only possible after the clock process $W$ is brought into the picture.
Recall the sets $\sD\subset \sD^*\subset D$ defined in (\ref{D}).

Next, we denote $\hC^{1,2,1}_0(D)$ to be the set of all functions $\vf\in \hC^{1,2,1}(\sD)$ such that
for $\eta=\vf$, $\vf_t$, $\vf_x$, $\vf_{xx}$, $\vf_w$, it holds that $\lim_{(t,y,v)\to (s,x,w)\atop  (t,y,v)\in\sD}\eta(t,y,v)=\eta(s,x,w)$,
for all $(s,x,w)\in D$;
and $\vf(s,x,w)=0$, for $(s,x,w)\notin D$.  We note that while
a function $\vf\in\hC^{1,2,1}_0(D)$ is well-defined on $D$, it is not necessarily continuous on
the boundaries $\{(s,x,w): x=0 ~\mbox{or}~ w=0 ~\mbox{or}~ w=s\}$.

Next, we define the following function:
\bea
\label{H0}
H(s,x,w,u, \xi, A, z, \g, a):=\frac{\sigma^{2}}2\gamma^{2} x^{2}A+(p+r x-a)\xi^1
			+\xi^2+ \l(w)z+(a-cu),
\eea
where $\xi=(\xi^1,\xi^2)\in \hR^2$, $u, A, z\in\hR$, and $(\g,a)\in [0,1]\times [0,M]$. For
 $\vf\in\hC^{1,2,1}_0(D)$, we define the following Hamiltonian:
\bea
\label{H}
\sH(s,x,w,\vf,\vf_x, \vf_w, \vf_{xx},\g, a):=H(s,x,w, \vf, \nabla \vf, \vf_{xx}, I(\vf), \g, a),
\eea
where $\nabla\vf:=(\vf_x,\vf_w)$ and  $I[\vf]$ is the integral operator defined by
\bea
\label{Iv}
 I[\vf]:=\int_0^\infty[\vf(s,x-u,0)-\vf(s,x,w)]dG(u)=\int_0^x\vf(s,x-u,0)dG(u)-\vf(s,x,w).
 \eea
Here the last 	equality is due to the fact that $\vf(s,x,w)=0$ for $x<0$.
The main purpose of this section is to show that the value function $V$ is a viscosity solution of the following
HJB equation:
		\bea
			\label{HJB}
			\left\{\ba{lll}
			\{V_s+\sL[V]\}(s,x,w)=0; \q (s,x,w)\in \sD;\\
			V(T, x, w)=0,
			\ea\right.
			\eea
where $\sL[\,\cd\,]$ is the second-order  partial integro-differential operator: for $\vf\in\hC^{1,2,1}_0(D)$,
\bea
\label{sL}
\sL[\vf](s,x,w):=\sup_{\gamma\in [0,1], a\in [0,M]}\sH(s,x,w,\vf,\vf_x,\vf_w, \vf_{xx}, \g, a).
\eea
\begin{rem}
\label{C121}
{\rm
(i) As we pointed out before, even a classical solution to the HJB equation (\ref{HJB}) may have discontinuity on the
boundaries $\{x=0\}$ or $\{w=0\}$ or $\{w=s\}$, and (\ref{HJB}) only specifies the boundary value of $\sD$ at
$s=T$.

(ii) To guarantee the well-posedness we shall consider the {\it constrained} viscosity solutions (cf. e.g.,
\cite{Soner}), for which the following observation is crucial.
Let $V\in \hC^{1,2,1}_0(D)$ be a classical solution so that (\ref{HJB}) holds on $\sD^*$. Consider the  point $(s,0,w)\in \pa \sD^*$. Let
$\vf\in\hC^{1,2,1}_0(D)$ be such that $0=[V-\vf](s,0,w)=\max_{(t,y,v)\in \sD^*}[V-\vf](t,y,v)$. Then one must have
$(\pa_t, \nabla)(V-\vf)(s,0,w)=\alpha\nu$ for some $\a>0$, where $\nabla=(\pa_x,\pa_w)$ and $\n$ is the outward normal vector of $\sD^*$ at the boundary $\{x=0\}$ (i.e.,
$\nu=(0,-1,0)$), and $I[V-\vf](s,0,w)=-[V-\vf](s,0,w)=0$ since $[V-\vf](s,y,w)=0$ for $y\le 0$. Thus,  for any $(\g, a)\in[0,1]\times [0,M]$ we obtain that
\bea
\label{vfleV0}
&&[\vf_s+\sH(\cd,\vf,\vf_x, \vf_w,\vf_{xx},\gamma,a)](s,0,w)
=\Big[\vf_s+((p-a,1), \nabla\vf)+ \lambda I[\vf]+(a-c\vf)\Big](s,0,w)\nonumber\\
&& \qq\qq\qq\qq\qq=[V_s+\sH(\cd,V,\nabla V, V_{xx}, I(V), \gamma,a)](s,0,w)+\a(p-a).
\eea
 Consequently, assuming $a\le p$ (which is natural in the case $x=0$!) we have
 \bea
 \label{vfleV}
\{ \vf_s+\sL[\vf]\}(s,0,w)\ge \{ V_s+\sL[V]\}(s,0,w)=0,
\eea
For the other two boundaries $\{w=0\}$ and $\{w=s\}$, we note that  $[V_{xx}-\vf_{xx}]\le 0$ and the corresponding outward normal vectors are
$\n=(0,0,-1)$ and $(-1,0,1)$, respectively. Therefore, a similar calculation as (\ref{vfleV0}), noting that $((1,p+rx-a,1), \n)=-1, 0$, respectively, would lead to (\ref{vfleV}) in both cases. In other words, we can extend the ``subsolution property" of (\ref{HJB}) to $\sD^*$.
%
\qed}
\end{rem}

We are now ready to give the definition of the so-called {\it constrained} viscosity solution.
	\begin{defn}
		Let ${\cal O}\subseteq \sD^*$ be a subset such that  $\pa_T\cO:=\{(T,y,v)\in\pa\cO\}\neq \emptyset$, and let $v\in \hC(\cO)$.
		
		(a) We say that $v$ is a viscosity subsolution of (\ref{HJB}) on ${\cal O}$,   if $v(T,y,v)\le 0$, for $(T, y, v)\in \pa_T\cO$;
 and for any $(s,x,w)\in \cO$ and $\vf\in \hC^{1,2,1}_0(\cO)$
such that $0=[v-\vf](s,x,w)=\max_{(t,y,v)\in\cO}[v-\vf](t,y,v)$, it holds that
\bea
\label{vsub}	
			\vf_s(s,x,w)+\sL[\vf](s,x,w)\geq 0.
\eea

		(b) We say that $v$ is a viscosity supersolution of (\ref{HJB}) on ${\cal O}$,   if $v(T,y,v)\ge 0$, for $(T, y, v)\in \pa_T\cO$;
 and for any $(s,x,w)\in\cO  $ and $\vf\in \hC^{1,2,1}_0(\cO)$ such that $0=[v-\p](s,x,w)=\min_{(t,y,v)\in\cO}[v-\vf](t,y,v)$,
it holds that
\bea
\label{vsup}	
			\vf_s(s,x,w)+\sL[\vf](s,x,w)\leq 0.
\eea
		
		(c) We say that $v\in \hC(D)$ is a ``constrained viscosity solution" of (\ref{HJB}) on $\sD^*$ if it is both a  viscosity subsolution of (\ref{HJB}) on $\sD^*$ and a viscosity supersolution of (\ref{HJB}) on ${\sD}$.
\qed
	\end{defn}
	
\begin{rem}
\label{vis}
{\rm (i) We note that the main feature of the constrained viscosity solution is that its subsolution is defined  on
	$\sD^*$, which is justified in Remark \ref{C121}-(ii). This turns out to be  essential for the
	comparison theorem, whence the uniqueness.
	
	(ii) The inequalities in (\ref{vsub}) and (\ref{vsup}) are opposite than the usual sub- and super-solutions, due to the fact that the HJB equation (\ref{HJB}) 	is a terminal value problem.}
\qed
	\end{rem}
		
		As in the viscosity theory, it is often convenient to define viscosity solution in terms of the sub-(super-) differentials, (or parabolic 	 sub-(super-)jets). To this end we introduce the following notions:
	\begin{defn}
	\label{jets}		
	Let $\cO\subseteq \sD^*$, $u\in \hC({\cal O})$, and $(s,x,w)\in {\cal O}$. The set of parabolic super-jets of $u$
		at $(s,x,w)$, denoted by ${\sP}_{{\cal O}}^{+(1,2,1)}u(s,x,w)$,  is defined as the set of all 	$(q,\xi,A)\in \hR\times \hR^{2}\times \hR$ such that for all 	 $(s,X):=(s, x,w), (t,Y):=(t,y,v) \in\cO$, it holds that
	\bea
	\label{supjet}
	u(t,Y)\leq u(s,X)+q(t-s)+(\xi, Y-X)+\frac{1}{2}A(x-y)^2+o(|t-s|+|w-v|+|y-x|^{2}),
	\eea


The set of parabolic sub-jets of $u$ at $(s,x,w)\in\cO$, denoted by ${\sP}_{{\cal O}}^{-(1,2,1)}u(s,x,w)$, is the set of all
$(q,p,A)\in \hR\times \hR^{2}\times \hR$  such that (\ref{supjet}) holds with ``$\,\le$" being replaced by ``$\,\ge$".
\qed		
	\end{defn}
The closure of  ${{\sP}}_{{\cal O}}^{+(1,2,1)}u(s,x,w)$ (resp. ${\sP}_{{\cal O}}^{-(1,2,1)}u(s,x,w)$), denoted by
$\bar{{\sP}}_{{\cal O}}^{+(1,2,1)}u(s,x,w)$ (resp. $\bar{{\sP}}_{{\cal O}}^{-(1,2,1)}u(s,x,w)$), is defined as the set of
all $(q,\xi,A)\in \hR\times \hR^{2}\times \hR$ such that there exists
$(s_{n},x_{n},w_{n})\in {\cal O}$ and $(q_{n},\xi_{n},A_{n})\in {{\sP}}_{{\cal O}}^{+(1,2,1)}u(s_{n},x_{n},w_{n})$ (resp.
${{\sP}}_{{\cal O}}^{-(1,2,1)}u(s_{n},x_{n},w_{n})$), and that
$ ((s_{n},x_{n},w_{n}),u(s_{n},x_{n},w_{n}),q_{n},\xi_{n},A_{n})\rightarrow ((s,x,w),u(s,x,w),q,\xi,A)$, as $n\to\infty$.

	We now define the constrained viscosity solution in terms of the parabolic jets. The equivalence between the two definitions
	in such a setting can be found in, for example, \cite{AlT, Sayah}.
	\begin{defn}
	\label{VD2}
	Let $\cO\subseteq \sD^*$, $u\in \hC({\cal O})$. We say that $\underline{u}$  (resp. $\bar{u}\in \hC(\cO)$ is a viscosity subsolution (resp. supersolution) of (\ref{HJB})
	on $\cO$ if for any $(s,x,w)\in\cO$, it holds that
		\beaa
	\label{sub}
		q+\sup_{\gamma\in [0,1], a\in [0,M]}H(s,x,w, \underline{u},\xi,A, I[\underline{u}], \gamma,a)&\geq &0 \\
		(\mbox{resp.}\qq  q+\sup_{\gamma\in [0,1], a\in [0,M]}H(s,x,w, \bar{u},\xi,A, I[\bar{u}], \gamma,a)&\leq &0),		
		\eeaa
for all $(q,(p^1,p^2),A)\in {\sP}_{\cO}^{+(1,2,1)}\underline{u}(s,x,w)$ (resp. ${\sP}_{{{\cal O}}}^{-(1,2,1)}\bar{u}(s,x,w)$).

In particular, we say that $u$ is a ``constrained  viscosity solution" of (\ref{HJB}) on $\sD^*$ if it is both a viscosity subsolution on $\sD^*$, and a
viscosity supersolution on $\sD$.
\qed
	\end{defn}
	
	In the rest of the paper, we shall assume that  all solutions of (\ref{HJB}) satisfy
$u(s,x,w)=0$,  for $(s,x,w)\notin D$. We now give the main result of this section.
	\begin{thm}
		Assume that Assumption \ref{assump1} is in force. Then, the value function $V$ of problem (\ref{V1}) is a constrained  viscosity solution of (\ref{HJB}) on $\sD^*$.
	\end{thm}
	
	{\it Proof.}  {\it  Supersolution.} Given $(s,x,w)\in\sD$. Let $\varphi \in \hC^{1,2,1}_0(D)$ be such that $V-\varphi$ attains its minimum  at $(s,x,w)$ with $\varphi(s,x,w)=V(s,x,w)$. For any $h>0$ such that
	$s\le s+h<T$, let us denote $\tau^h_s:=s+h\wedge T^{s,w}_1$, 	and $\tilde U_{1}=\D Q^{s,w}_{T^{s,w}_1}$.
For any $(\g_0,a_{0})\in[0,1]\times[0,M]$, we consider the following ``feedback" strategy:
$\pi^0_t=(\g_0, a_0{\bf 1}_{\{t<\t_0\}}+p{\bf 1}_{\{t\ge \t_0\}})$, $t\in[s,T]$, where $\t_0=\inf\{t>s, X^{\pi^0}_t=0\}$.
Then $\pi^0 \in \sU_{ad}[s,T]$,  and it is readily seen from (\ref{Xsol}) that ruin can only happen at a jump time, that is,
$T^{s,w}_1\le \t^{\pi^0}_s$, and $R_{t}^{0}:=(t, X^{\pi^0,s,w,x}_t, W^{s,w}_t)\in\sD$, for $t\in[s, \t^h_s)$.
	
	Next, by DPP (Theorem \ref{DPPthm})  and the properties of $\vf$ we have
\bea
\label{supest}
		0&\geq & \hE_{sxw}\Big[\int_{s}^{\tau^h_s}e^{-c(t-s)} (a_0{\bf 1}_{\{t<\t_0\}}+p{\bf 1}_{\{t\ge \t_0\}})dt+e^{-c(\tau^h_s -s)}V(R^0_{\tau^h_s})  \Big]-V(s,x,w)\nonumber\\
				 &\geq&\hE_{sxw}\Big[\int_{s}^{\tau^h_s}e^{-c(t-s)}a_{0}dt{\bf 1}_{\{\t^h_s<\t_0\}}+e^{-c(\tau^h_s -s)}\varphi(R^{0}_{\tau^h_s}) \Big)\Big]-\varphi(s,x,w)\\
		&=&\hE_{sxw}\Big[\int_{s}^{\tau^h_s}e^{-c(t-s)}a_{0}dt{\bf 1}_{\{\t^h_s<\t_0\}}\Big]+\hE_{sxw}\Big[e^{-c(\t^h_s -s)}[\vf(R^{0}_{\tau^h_s})-\vf(R^{0}_{\tau^h_s-})]
		{\bf 1}_{\{T^{s,w}_1<h\}}\Big]\nonumber\\
		&&+\hE_{sxw}\Big[e^{-c(\tau^h_s -s)}\vf(R^{0}_{\tau^h_s-})-\vf(s,x,w)\Big]:= I_1+I_2+I_3,\nonumber
		\eea
where $I_i$, $i=1,2,3$ are the three terms on the right hand side above. Clearly,  we have
\bea
	\label{supI1}
		I_1=\frac{a_0}c \Big\{[1-e^{-ch}]\hP(\t_0>s+h, T^{s,w}_1>h)+ \int_0^h\Big[1-e^{-ct}]\hP\{\t_0>s+t)dF_{T^{s,w}_1}(t)\Big\},
	\eea
Since $\t^h_s=s+T^{s,w}_1$ on $\{T^{s,w}_1<h\}$, we have
\bea
\label{supI2}
 I_2&=& \hE_{sxw}\Big[e^{-cT^{s,w}_1}[\vf(R^0_{s+T^{s,w}_1})-\vf(R^0_{s+T^{s,w}_1)-})]{\bf 1}_{\{T^{s,w}_1<h\}}\Big]\\
 &=&\hE_{sxw}\Big[\int_0^\infty\neg\neg\int_0^{h} e^{-ct}[\vf(s+t,X^{\pi^0}_{(s+t)-}-u, 0)-\vf(t, X^{\pi^0}_{(s+t)-}, W^{s,w}_{(s+t)-})]dF_{T^{s,w}_1}(t)dG(u)\Big].\nonumber
  \eea
 Since there is no jumps on $[s, \t^h_s)$, applying It\^o's formula (and denoting $\th(x):=rx+p$) we get
  \bea
  \label{supI3}
        I_3& =&\hE_{sxw}\Big[\int_{s}^{\tau^h_s}e^{-c(t-s)}\big\{\neg-\neg c\vf+ \vf_t+ ((\th(X^{\pi^0}_{t})-a_{0}, 1), \nabla\vf)
        +\frac{(\si\g_0X^{\pi^0}_t)^2}{2}\varphi_{xx} ^{2}\big\}(R^0_{t})dt\Big]\\
		&=& \hE_{sxw}\Big[\int_{0}^{h}\neg\neg{\bf 1}_{\{T^{s,w}_1\ge t\}}e^{-ct}\big\{\neg-\neg c\vf+ \vf_t+
		((\th(X^{\pi^0}_{s+t})-a_{0}, 1), \nabla\vf)+\frac{(\si\g_0X^{\pi^0}_{s+t})^2}{2}\varphi_{xx} ^{2}\big\}(R^0_{s+t})dt\Big]
		\nonumber \\
		&=&\hE_{sxw}\Big[\int_{0}^{h}\bar F_{T^{s,w}_1}(t)e^{-ct}\big\{\neg-\neg c\vf+
		 \vf_t+ ((\th(X^{\pi^0}_{s+t})-a_{0}, 1), \nabla\vf)+\frac{(\si\g_0X^{\pi^0}_{s+t})^2}{2}\varphi_{xx} ^{2}\big\}(R^0_{s+t})dt\Big].
		\nonumber 		
	\eea
Recall that $dF_{T^{s,w}_1}(t)=\l(w)\bar F_{T^{s,w}_1}(t)dt=\l(w)e^{-\int_w^{w+t}\l(u)du}dt$, and
$\bar F_{T^{s,w}_1}(0)=1$, dividing both sides of (\ref{supest}) by $h$ and then sending $h$ to $0$ we obtain, in light of
(\ref{supI1})--(\ref{supI3}),
	\bea
	\label{supsol}
	0\ge \{\vf_t +\sH(\cd, \vf, \vf_x,\vf_w,\vf_{xx}, \g_0, a_0)\}(s,x,w).
	\eea
Since $(\g_0,a_0)$ is arbitrary, we conclude that  $V$ is a viscosity supersolution on $\sD$.

\ms
{\it Subsolution.}  We shall now argue that $V$ is a viscosity subsolution on  $\sD^*$. Suppose not, then we shall first show that
there exist $(s,x,w)\in \sD^*$, $\p\in \hC^{1,2,1}_0(D)$, and constants $\e>0$, $\rho>0$, such that $0=[V-\p](s,x,w)=\max_{(t,y,v)\in\sD^*} [V-\p](t,y,v)$, but
\bea
	\label{Vlep}	
	\left\{\ba{lll}
		\{\p_s+\sL[\p]\}(t,y,v) \leq  -\e c,  \qq & (t,y,v)\in \overline{B_{\rho}(s,x,w)\cap {\sD^*}}\setminus \{t=T\};\\
V(t,y,v)\leq \psi(t,y,v)-\e, \qq & (t,y,v)\in \partial B_{\rho}(s,x,w)\cap \sD^*,
\ea\right.
\eea
where $B_{\rho}(s,x,w)$ is the  open ball  centered at  $(s,x,w)$ with radius $\rho$. 
To see this, we note that if $V$ is not a viscosity subsolution on $\sD^*$, then there must exist $(s,x,w)\in\sD^*$ and
$\p^0\in \hC^{1,2,1}_0(D)$, such that $0=[V-\p^0](s,x,w)=\max_{(t,y,v)\in\sD^*} [V-\p^0](t,y,v)$, but
\bea
\label{Llee0}
 \{\p^0_s+\sL[\p^0]\}(s,x,w) =-2\eta<0, \qq \mbox{for some $\eta>0$.}
 \eea	
We shall consider two cases.

{\it Case 1. $x>0$.} In this case we introduce the function:
	\bea
\label{peta}
\psi(t,y,v):=\psi^{0}(t,y,v)+\frac{ \eta [(t-s)^{2}+(y-x)^{2}+(v-w)^{2}]^{2}}{\l(w)(x^{2}+w^{2})^{2}}, \qq  (t,y,v)\in D.
\eea
Clearly,  $\psi\in \hC^{1,2,1}_0(D)$,
$\psi(s,x,w)=\psi^{0}(s,x,w)=V(s,x,w)$,   and $\p(t,y,v)>V(t,y,v)$, for all  $(t,y,v)\in D\setminus (s,x,w)$. Furthermore, it is easy to check that $(\p_s, \nabla\p)(s,x,w)=(\p^0_s, \nabla\p^0)(s,x,w)$, $\p_{yy}(s,x,w)=\p^0_{yy}(s,x,w)$,
and
\beaa
	 \label{61602}
		 \l(w)\int_{0}^{x} \psi(s,x-u,0)dG(u)
		\le \l(w)\int_{0}^{x} \psi^{0}(s,x-u,0)dG(u)+\eta.
	\eeaa
Consequently, we see that
$$\{\p_s+\sL[\p]\}(s,x,w)\le \{\p^0_s+\sL[\p^0]\}(s,x,w)+\eta=-\eta<0.
$$
By continuity of $\p_s+\sL[\p]$, we can then find $\rho>0$ such that
\bea
\label{Llee}
		\{\p_{t}+\sL[\p]\}(t,y,v)<-\eta/ 2, \qq \mbox{for $(t,y,v)\in\overline{B_{\rho}(s,x,w)\cap {\sD^*}}\setminus \{t=T\}$.}
\eea
Note also that for $(t,y,v)\in \partial B_{\rho}(s,x,w)\cap \sD^*$, one has
	\bea
	\label{Vlep1}
		V(t,y,v)&\leq &
		\psi(t,y,v)-\frac{ \eta \rho^{4}}{\l(w)(x^{2}+w^{2})^{2}}.
	\eea
Thus if we choose $\e={\rm min}\Big\{\frac{\eta}{2c}, \frac{ \eta \rho^{4}}{\lambda(w)(x^{2}+w^{2})^{2}}\Big\}$,  then (\ref{Llee}) and (\ref{Vlep1})
become (\ref{Vlep}).

\ms
{\it Case 2. $x=0$.} In this case we have
	$$	\p^0_s-\sL[\p^0](s,0,w)= \sup_{a\in [0,M]}[((1,p-a,1), (\p^0_s, \nabla\p^0))(s,0,w)- (c+\l(w))\p^0(s,0,w)+a].
	$$
If we define $\p(t,y,v)= \p^{0}(t,y,v)+\eta[ (t-s)^{2}+y^{2}+(v-w)^{2}]$, for $(t, y, v)\in D$, and $\e={\rm min}\Big\{\frac{\eta}{2c},  \rho^{2}\Big\}$,
then a similar calculation as before shows that (\ref{Vlep}) still holds,
proving the claim. In what follows we shall argue that this will lead to a contradiction.

To this end, fix any $\pi=(\gamma,a)\in \sU_{ad}^{s,w}[s,T]$, and let $R^{s,x,w}_t=(t,X^{s,x,w}_{t},W^{s,w}_{t})$. Define $\tau_{\rho}:=\inf\{t>s: R_t\notin \overline{B_{\rho}(s,x,w)\cap {\sD^*}}\}$,  $\tau:=\tau_{\rho}\wedge T^{s,w}_{1}$, and
denote 
$R_{t}=R^{s,x,w}_t$ for simplicity. Applying It\^o's formula to $e^{-c(t-s)}\p(R_t)$ from $s$ to $\t$
we have
 	\bea
	\label{VpIto}
			&& \int_{s}^{\tau}e^{-c(t-s)}a_{t}dt+ e^{-c(\tau -s )} V(R_{\tau })=\int_{s}^{\tau}e^{-c(t-s)}a_{t}dt+ e^{-c(\tau -s )} [\p(R_{\tau })+(V(R_{\tau })-\p(R_\t))]\nonumber\\
		&=& e^{-c(\tau -s )} [V(R_{\tau })-\p(R_\t)]+\psi(s,x,w)\nonumber\\
	&&	+ \int_{s}^{\tau}e^{-c(t-s)}[a_{t}-c \psi+\psi_{t}+\psi_{w}+(rX_{t}+p-a_{t})\psi_{x}
		+\frac{1}{2} X_{t}^{2}\sigma^{2}\gamma_{t}^{2}\psi_{xx}](R_{t})dt \\
		&&		+\int_{s}^{\tau}e^{-c(t-s)}\psi_{x}(R_{t})\sigma \gamma_{t}X_{t}dW_{t}  +\sum_{s\leq t\leq \tau}e^{-c(t-s)}(\psi(R_t)-\psi(R_{t-})). \nonumber
	\end{eqnarray}
Then, on the set $\{\tau_{\rho}\geq T^{s,w}_{1}\}$, we have $\tau= T^{s,w}_{1}$. Since the ruin only happens at the claim arrival
times, we have $\t^\pi\ge T^{s,w}_1$. In the case that $\tau^{\pi}= T^{s,w}_1$,  $X_{ T^{s,w}_1}<0$ and  $V(R_{T^{s,w}_1})=\psi(R_{T^{s,w}_1})=0$; whereas in the case  $\tau^{\pi}>T^{s,w}_{1}$, we have $R_{T^{s,w}_1}\in D$, and $V(R_{T^{s,w}_1})\leq\psi(R_{T^{s,w}_1})$. 

On the other hand, we note that on the set $\{\tau_{\rho}<T^{s,w}_{1}\}$, $\t=\t_\rho$, and since $(\tau_{\rho},X_{\tau_{\rho}},W_{\tau_{\rho}})\in \partial B_{\rho}(s,x,w)\cap {{\sD^*}}$, we derive from (\ref{Vlep}) that  $[V(R_{\tau_{\rho} })-\p(R_{\tau_{\rho}})]\leq -\varepsilon $.   Thus, noting that $W_{T^{s,w}_1}=0$, and that both $\p_x$ and $\g$ are bounded,
 we deduce from (\ref{VpIto}) that
 \begin{eqnarray}
	\label{subcontra}
		&&\hE_{sxw}\Big[\int_{s}^{\tau}e^{-c(t-s)}a_{t}dt+ e^{-c(\tau -s )} V(\tau , X_{\tau },W_{\tau })\Big] \nonumber \\
		&\leq& \hE\Big[\psi(s,x,w)-\e e^{-c(\tau_{\rho} -s )}{\bf 1}_{\{\tau_{\rho}<T^{s,w}_{1}\}}+\int_{s}^{\tau}e^{-c(t-s)} [\p_t+
		\sH(\cds, \g_t,a_t)(R_{t})]dt\Big]\\
		&\leq &\psi(s,x,w)-\e\hE_{sxw}\big[e^{-c(\tau -s )}{\bf 1}_{\{\tau_{\rho}<T^{s,w}_{1}\}}+(1-e^{-c(\tau-s)})\big]\nonumber\\
		&=&V(s,x,w) -\e\hE_{sxw}[(1-e^{-c(T^{s,w}_{1} -s)}{\bf 1}_{\{\tau_{\rho}\geq T^{s,w}_{1}\}})]
		\le V(s,x,w)-\e \hE_{sxw}(1-e^{-c( T^{s,w}_{1}-s)}).\nonumber
	\end{eqnarray}
Since $\hP\{T^{s,w}_1>s\}=1$, we see that (\ref{subcontra}) contradicts the Dynamic Programming Principle (\ref{DPP}).
This completes the proof.
\qed 

\section{Comparison Principle and Uniqueness}
\setcounter{equation}{0}

In this section, we present a Comparison Theorem that would imply the uniqueness among a certain class of the constrained viscosity solutions of (\ref{HJB}) to which the value function belong. To be more precise, we introduce to following subset of $\hC(D)$.
 \begin{defn}
\label{cond}
We say that a function $u\in \hC(D)$ is of class (L) if

 (i) $u(s,x,w)\ge 0$, $(s,x,w)\in D$, and $u$ is  uniformly continuous on   $D$;

(ii)  the mapping $x\mapsto u(s,x,w)$ is increasing, and  $\lim_{x\rightarrow\infty}u(s,x,w)=\frac{M}{c}[1-e^{-c(T-s)}]$;

(iii)  $u(T,y,v)=0$ for any $(y, v)\in [0,\infty)\times[0,T]$.
\end{defn}

Clearly, the value function $V$ of problem (\ref{V1}) is of class (L), thanks to Proposition \ref{vb},  Proposition \ref{conts},  Theorem \ref{continuityx}, and
Corollary \ref{coninw}. Our goal is to show that following {\it Comparison Principle}.

\begin{thm}[Comparison Principle]
\label{comp}
		Assume that Assumption \ref{assump1} is in force. Let $\underline{u}$ be a viscosity subsolution of (\ref{HJB}) on $\sD^*$ and
		$\bar{u}$ be a viscosity supersolution of (\ref{HJB}) on $\sD$. If both $\bar{u}$ and $\underline{u}$ are of class (L), then
		$\underline{u}\leq \bar{u}$ on $D$.
	
	Consequently, there is at most one constrained viscosity solution of class (L) to (\ref{HJB})  on $D$.
	\end{thm}
%

		
{\it Proof.}    We first perturb the supersolution slightly so that all the inequalities involved become strict. Define, for $\rho>1$,  $\theta, \vsi>0$,
$$\bar{u}^{\rho,\theta,\vsi}(t,y, v)=\rho\bar{u}(t,y, v)+\theta\frac{T-t+\vsi}{t}.$$
Then it is straightforward to check that $\bar{u}^{\rho,\theta,\vsi}(t,y, v)$ is also a supersolution of (\ref{HJB}) on ${\sD}$. In fact, it is easy to see that $\rho\bar{u}$ is a supersolution of (\ref{HJB}) in ${\sD}$ as $\rho>1$, and for any $(s,x,w)\in \sD$ and $\varphi\in C^{1,2,1}_0(D)$ such that $0=[\bar{u}^{\rho, \theta, \vsi}-\vf](s,x,w)=\min_{(t,y,v)\in \sD}[\bar{u}^{\rho, \theta, \vsi}-\vf](t,y,v)$,
it holds that
	\begin{eqnarray*}
		[\vf_t+\sup_{\g, a}\sH(\cd, \bar{u}^{\rho,\theta, \vsi}, \vf_x,\vf_w,\vf_{xx}, \g, \a)](s,x,w)
		\leq [\vf_t+\sup_{\g,a}\sH(\cd, \rho\bar{u}, \tilde\vf_x,\tilde\vf_w,\tilde\vf_{xx}, \g, a)](s,x,w)\le 0,
	\end{eqnarray*}
where $\tilde\vf(t,y,v):=\varphi(t,y,v)-{\theta(T-t+\varsigma)}/{t}$, i.e.,  $\bar{u}^{\rho, \theta, \vsi}$ is a viscosity supersolution on $\sD$. We shall argue that $\underline{u}\leq\bar{u}^{\rho,\theta}$, which will lead to the desired
comparison result as $\lim_{\rho\da 0,\theta\da 0, \vsi\da 0}\bar{u}^{\rho,\theta,\vsi}=\bar u$.

To this end, we first note that
	$\lim_{t\rightarrow 0}	\bar{u}^{\rho,\theta}(t,y,v)=+\infty$. Thus we need only show that
 $\underline{u}\leq\bar{u}^{\rho,\theta}$ on $\sD^*\setminus \{t=0\}$. Next, note that both $\underline u$ and $\bar u $ are of
 class (L), we have  (recall Definition \ref{cond})
 	\bea
	\label{boundary1}
\lim_{y\rightarrow\infty}(\underline{u}(t,y,v)-\bar{u}^{\rho,\theta,\vsi}(t,y,v))=(1-\rho)\frac{M}{c}[1-e^{-c(T-t)}]-\frac{\theta(T-t+\varsigma)}{t}\le -\frac{\th \vsi}{T}<0,
	\eea
for all $0< t\le T$. Thus, by Dini's Theorem, the convergence in (\ref{boundary1}) is uniform in $(t, y, v)$, and we can choose $b>0$ so that $\underline{u}(t,y,v)<\bar{u}^{\rho,\theta}(t,y,v)$  for $y\ge b$, $0< t <T$,
and $0\leq v \leq t$. Consequently,  it suffices to show that
\bea
\label{D0}
\underline{u}(t,y,v)\leq\bar{u}^{\rho,\theta,\vsi}(t,y,v), \q \mbox{on}~
{\sD_b}=\{(t,y,v):0< t <T, 0\leq y < b, 0\leq v \leq t\}.
\eea

Suppose (\ref{D0}) is not true, then there exists $(t^{*}, y^{*}, v^{*})\in \bar{\sD_b}$ such that
	\bea
	\label{Mb}
		M_b:=\sup_{\sD_b}(\underline{u}(t,y,v)-\bar{u}^{\rho,\theta, \vsi}(t,y,v))=\underline{u}(t^{*}, y^{*}, v^{*})-\bar{u}^{\rho,\theta, \vsi}(t^{*}, y^{*}, v^{*})>0.
		\eea

Next, we denote $\sD^0_b:=$int$\sD_b$, and
\bea
\label{D02}
\sD^1_b :=\pa\sD_b\cap \sD_b =\pa \sD_b\setminus [\{t=0\}\cup \{t=T\}\cup \{y=b\}].
\eea
We note that $\underline{u}(t,y,v)-\bar{u}^{\rho, \theta,\vsi}(t,y,v)\le 0$, for $t=0, T$ or $y=b$, thus $(t^*, y^*, v^*)$ can only happen on
$\sD^0_b\cup \sD^1_b$. We shall consider the following two cases separately.

{\it Case 1.} We assume that $(t^{*}, y^{*}, v^{*})\in \sD^0_b$, but
\bea
\label{D1b}
\underline{u}(t,y,v)-\bar{u}^{\rho,\theta, \vsi}(t,y,v)<M_b, \qq (t,y,v)\in\sD^1_b.
\eea
In this case we follow a more or less standard argument. For  $\e>0$, we define an auxiliary function:
	\bea
	\label{Sigma0}
		\Sigma^b_{\varepsilon}(t,x,w,y,v)=\underline{u}(t,x,w)-\bar{u}^{\rho,\theta, \vsi}(t,y,v)-\frac{1}{2\varepsilon}(x-y)^{2}-\frac{1}{2 \varepsilon} (w-v)^{2},
	\eea
for $(t,x,w,y,v)\in\sC_b:=\{(t,x,w,y,v): t\in[0,T], x,y\in[0,b], w,v\in[0,t]\}$.
Since $\sC_b$ is compact, there exist  $\{(t_\e, x_\e,w_\e,y_\e,v_\e)\}_{\e>0}\subset\sC_b$, such that
\bea
\label{Me}
M_{\e,b}:=\max_{\sC_b}\Sigma^b_{\e}(t,x,w,y,v)=\Sigma^b_\e(t_\e, x_\e,w_\e, y_\e,v_\e).
\eea
We claim that for some $\e_0>0$,
$(t_\e, x_\e,w_\e,y_\e,v_\e))\in$ int$\,\sC_b$, whenever $0<\e<\e_0$.

Indeed, suppose not,
then there is a sequence  $\e_n\da 0$, such that $(t_{\e_{n}}, x_{\e_{n}},w_{\e_{n}},y_{\e_{n}},v_{\e_{n}})\in \pa\sC_b$,
the boundary of $\sC_b$, and that (\ref{Me}) holds for each $n$.
Now since $\pa\sC_b$ is compact, we can find a subsequence, may assume $(t_{\e_{n}}, x_{\e_{n}},w_{\e_{n}}, y_{\e_{n}},v_{\e_{n}}) $ itself, such that $(t_{\e_{n}}, x_{\e_{n}},w_{\e_{n}},y_{\e_{n}},v_{\e_{n}})\rightarrow (\hat{t},\hat{x},\hat{w},\hat{y},\hat{v})\in\pa\sC_b$.

Note that the function $\underline{u}$ is  continuous and bounded on $D$, and
\bea
\label{Mb1}
\Sigma^b_{\e_n}(t_{\e_{n}}, x_{\e_{n}},w_{\e_{n}},y_{\e_{n}},v_{\e_n})= M_{\e_n,b}\geq \Sigma^b_{\e_n}(t^{*}, y^{*}, v^{*};  y^{*}, v^{*})=M_b>0,
\eea
 it follows from (\ref{Sigma0}) and (\ref{Mb1}) that
	$$\frac{(x_{\e_{n}}-y_{\e_{n}})^{2}}{2\e_{n}} +\frac{(w_{\e_{n}}-v_{\e_{n}})^{2}}{2\e_{n}}\le  \underline{u}( t_{\e_{n}}, x_{\e_{n}},w_{\e_{n}})
 \le \frac{M}c.$$
Letting $n\to\infty$ we obtain that $\hat x=\hat y$, $\hat w=\hat v$, which implies,
by (\ref{Mb1}),
\bea
\label{Mb2}
\underline{u}(\hat t, \hat x, \hat w)-\bar{u}^{\rho,\theta, \vsi}(\hat t, \hat x, \hat w)=\Sigma^b_{\e}(\hat t, \hat x, \hat w, \hat x,\hat w) =\lim_{n\to\infty}
\Sigma^b_\e(t_{\e_{n}}, x_{\e_{n}},w_{\e_{n}},y_{\e_{n}},v_{\e_n})\ge M_b>0.
\eea
But as before we note that $\underline{u}(t,y,v)-\bar{u}^{\rho,\theta,\vsi}(t,y,v)\le 0$ for
$t=0$, $t=T$, and $y=b$,
we conclude that $\hat t\neq 0,T$, and $\hat x<b$. In other words, $(\hat t, \hat x,\hat w)\in \pa\sD^0_b\setminus (\{t=0\}\cup\{t=T\}\cup\{y=b\}]=\sD_b^1$.
This, together with (\ref{Mb2}), contradicts the assumption (\ref{D1b}).

In what follows we shall assume that $(t_{\e}, x_{\e},w_{\e},y_{\e},v_{\e})\in$ int$\,\sC_b$, $\forall\e>0$. Applying
\cite[Theorem 8.3]{CIL}
one shows that for any $\delta>0$, there exist $q=\hat{q}\in \hR$ and $A, B\in\cS^2$ such that
	\beaa
\left\{\ba{lll}
\dis (q, ((x_{\e}-y_{\e})/\e, (w_{\e}-v_{\e})/\e ), A)\in \bar{{\sP}}_{\sD^0_b}^{1,2,+}\underline{u}(t_{\e},x_{\e},w_{\e}), \\
\dis	(\hat q, ((x_{\e}-y_{\e})/\e, (w_{\e}-v_{\e})/\e ), B)\in \bar{{\sP}}_{\sD^0_b}^{1,2,-}\bar{u}^{\rho,\theta, \vsi}(t_{\e},y_{\e},v_{\e}),
\ea\right.
	\eeaa
where $\bar\sP^{1,2,+}_{\sD^0_b}u(t,x,w)$ and $\bar\sP^{1,2,^-}_{\sD^0_b} u(t,y,v)$
are the closures of the usual parabolic super-(sub-)jets
of the  function $u$ at $(t,x,w), (t,y,v)\in \sD^0_b$, respectively (see \cite{CIL}), such that
	\begin{eqnarray}
		\frac{1}{\varepsilon}\begin{pmatrix}
			I & -I\\
			-I & I
		\end{pmatrix}+\frac{2\delta}{\varepsilon^{2}}
		\begin{pmatrix}
			I & -I\\
			-I & I
		\end{pmatrix}
		\geq
		\begin{pmatrix}
			A & 0\\
			0 &-B
		\end{pmatrix}	
	\end{eqnarray}
	where $I$ is the $2\times 2$ identity matrix. 	Taking $\delta=\varepsilon$, we have
	\begin{eqnarray}
		\frac{3}{\varepsilon}\begin{pmatrix}
			I & -I\\
			-I & I
		\end{pmatrix}\geq
		\begin{pmatrix}
			A & 0\\
			0 & -B
		\end{pmatrix}.
		\label{91801}
	\end{eqnarray}
Note that if we denote $A=[A_{ij}]_{i,j=1}^2$ and $B=[B_{ij}]_{i,j=1}^2$ and $\xi_\e:=( (x_{\e}-y_{\e})/\e, (w_{\e}-v_{\e})/\e)$, then $(q,\xi_\e, A)\in \bar\sP^{1,2,+}_{\sD^0_b}\underline{u}(t_\e,x_\e,w_\e)$,
(resp. $(\hat q, \xi_\e, B)\in \bar\sP^{1,2,-}_{\sD^0_b} \bar{u}^{\rho,\theta,\vsi}(t_\e,y_\e,v_\e)$) implies that $(q,\xi_\e, A_{11})\in \bar\sP^{+(1,2,1)}_{\sD^*}\underline{u}(t_\e,x_\e,w_\e)$ (resp. $(\hat q, \xi_\e, B_{11})\in \bar\sP^{-(1,2,1)}_{\sD} \bar{u}^{\rho,\theta,\vsi}(t_\e,y_\e,v_\e)$).
Since the functions $\underline{u}$, $\bar{u}^{\rho, \th, \vsi}$, and $H$ are all continuous in all variables, we may assume without
loss of generality that $(q,\xi_\e, A_{11})\in \sP^{+(1,2,1)}_{\sD^*}\underline{u}(t_\e,x_\e,w_\e)$ (resp. $(\hat q, \xi_\e, B_{11})\in \sP^{-(1,2,1)}_{\sD} \bar{u}^{\rho,\theta,\vsi}(t_\e,y_\e,v_\e)$) and, by Definition \ref{VD2},
	$$
\left\{\ba{lll}
\dis q+ \sup_{\gamma\in [0,1], a\in [0,M]}H(t_{\e},x_{\e},w_{\e},\underline{u},\xi_{\e}, A_{11}, I[\underline{u}], \g, a)\geq 0,
\\
\dis q+\sup_{\gamma\in [0,1], a\in [0,M]}H(t_{\e},y_{\e},v_{\e},\bar{u}^{\rho,\theta,\vsi},\xi_{\e}, B_{11}, I[\bar{u}^{\rho,\th,\vsi}], \g, a)\leq0.
\ea\right.
	$$
Furthermore, we note that
(\ref{91801}) in particular implies that
	\begin{eqnarray}
		A_{11}x_{\e}^{2}-B_{11}y_{\e}^{2}\leq \frac{3}{\e}(x_{\e}-y_{\e})^{2}.
	\end{eqnarray}
Thus, if we choose $(\gamma_{\e},a_{\e})\in {\rm arg max}_{(\gamma,a)\in [0,1]\times [0,M]}H(t_{\e},y_{\e},v_{\e},\underline{u},\xi_{\e}, A_{11}, I[\underline{u}],   \g, a)$, then we have
	\begin{eqnarray*}
	H(t_{\e},x_{\e},w_{\e},\underline{u},\xi_{\e}, A_{11}, \g_\e,a_\e)
	 -H(t_{\e}, y_\e,v_\e,\bar{u}^{\rho,\th,\vsi}, \xi_{\e}, B_{11}, \g_\e,a_\e)\geq 0.
	\end{eqnarray*}
Therefore, by definition (\ref{H}) we can easily deduce that
	\begin{eqnarray}
	\label{91802}
		&&c(\underline{u}(t_{\e},x_{\e},w_{\e})-\bar{u}^{\rho,\th,\vsi}(t_{\e}, y_{\e},v_{\e}))+\l(w_{\e})\underline{u}(t_{\e},x_{\e},w_{\e})-
		\l(v_{\e})\bar{u}^{\rho,\th,\vsi}(t_{\e}, y_{\e},v_{\e})\nonumber\\
		&\leq &\frac{1}{2}\sigma^{2}{\gamma_{\varepsilon}}^{2}(A_{11}x_{\e}^{2}-B_{11}y_{\e}^{2})
		+r\frac{(x_{\e}-y_{\e})^{2}}{\e} \nonumber \\
		&&+\l(w_{\e})\int_0^{x_{\e}} \underline{u}(t_{\e},x_{\e}-u,0)dG(u)-\l(v_{\e})\int_{0}^{y_{\e}} \bar{u}^{\rho,\theta,\vsi}(t_{\e},y_{\e}-u,0)dG(u)\\
		&\leq &\Big(\frac{3\sigma^{2}}2 	+r\Big)\frac{(x_{\e}-y_{\e})^{2}}{\e} \nonumber \\
		&&+\l(w_{\e})\int_{0}^{x_{\e}} \underline{u}(t_{\e},x_{\e}-u,0)dG(u)-\l(v_{\e})\int_{0}^{y_{\e}} \bar{u}^{\rho,\theta,\vsi}(t_{\e},y_{\e}-u,0)dG(u)\nonumber
	\end{eqnarray}
Now, again, since $(t_\e,x_\e,w_\e,y_\e,v_\e)\in \sC_b\subset \bar{\sC_b}$ which is compact, there exists a sequence
$\e_{m}\rightarrow 0$ such that
	$(t_{\e_{m}}, x_{\e_{m}},w_{\e_{m}}, y_{\e_{m}},v_{\e_{m}})\rightarrow (\bar{t},\bar{x},\bar{w},\bar{y},\bar{v})\in \bar{\sC_b}$.
By repeating the arguments before one shows that $\bar t\in(0,T)$, $\bar x=\bar y\in[0,b)$,  $\bar w=\bar v\in [0,t]$, i.e.,
and
 $$	
		\underline{u}(\bar{t},\bar{x},\bar{w})-\bar{u}^{\rho,\th,\vsi}(\bar{t}, \bar{x},\bar{w})=\lim_{\e_{m}\rightarrow 0}M_{\e_{m},b}\geq M_b,
$$
	we obtain that $(\bar t,\bar x,\bar w)\in \sD^0_b$. But on the other hand, replacing $\e$ by $\e_{m}$ and letting ${m}\to\infty$ in
(\ref{91802}) we have
	\begin{eqnarray*}
		(c+\l(\bar{w}))M_b
		&\leq & \l(\bar{w})\int_{0}^{\bar{x}}[\underline{u}(\bar{t},\bar{x}-u,0)- \bar{u}^{\rho,\th,\vsi}(\bar{t},\bar{x}-u,0)]dG(u)
		\le \l(\bar{w}) M_b.
	\end{eqnarray*}
	This is a contradiction as $c>0$ and $M_b>0$.

\ms
{\it Case 2.} We now consider the case $(t^{*}, y^{*}, v^{*})\in {\sD^1_b}$. We shall first move the point away from the boundary
$\sD^1_b$ into the interior $\sD^0_b$
and then argue as Case 1. To this end we borrow some arguments from
\cite{BKR}, \cite{Ishii} and \cite{Soner}.
First, since $(t^*, y^*, v^*)$ is on the boundary of  a simple polyhedron and $0<t^*<T$,   it is not hard to see that there exist $\eta=(\eta_{1},\eta_{2})\in \hR^{2}$, and $a>0$ such that
for any $(t, x, w)\in B^3_{a}(t^{*}, y^{*}, v^{*})\cap \sD^0_b $,  $0<\d \leq 1$, it holds that
	\bea
	\label{91601}	
	(t,y,v)\subset {\sD^0_b}, \qq {\rm whenever}\  (y,v)\in B^2_{\d a}(x+\d\eta_{1}, w+\d\eta_{2}).
	\eea
Here $B^n_{\rho}(\xi)$ denotes the ball centered at $\xi\in\hR^n$ with radius $\rho$.
	For any $\varepsilon>0$ and $0<\beta<1$, define the auxiliary functions: for $(t,x,w,y,v)\in \sC_b$,
	\begin{eqnarray*}
		\phi_{\e,\beta}(t,x,w,y,v):=\Big(\frac{x-y}{\sqrt{2\e}}+\beta\eta_{1}\Big)^{2}
		+\Big(\frac{w-v}{\sqrt{2\varepsilon}}+\beta\eta_{2}\Big)^{2}
		+\beta[(t-t^{*})^{2}+(x-y^{*})^{2}+ (w-v^{*})^{2}].
	\end{eqnarray*}
$\Sigma_{\e,\beta}(t,x,w,y,v):=\underline{u}(t,x,w)-\bar{u}^{\rho,\theta,\vsi}(t,y,v)-\phi_{\e, \beta}(t,x,w,y,v)$. Again, we have
\bea
\label{Meb}
M_{\e, \beta,b}:=\sup_{\sC_b}\Sigma_{\e,\beta}(t,x,w,y,v)\geq \Sigma_{\e,\beta}( t^{*}, y^{*}, v^{*}, y^{*}, v^{*})=M_b-\beta^{2}|\eta|^2>0,
\eea
for any $\e>0$ and $\beta<\beta_{0}$, for some $\beta_{0}>0$.
Now we fix $\beta\in(0,\beta_0)$ and denote, for simplicity, $( t_{\e}, x_{\e},w_{\e}, y_{\e},v_{\e})\in {\rm argmax}_{\sC_b} \Sigma_{\e, \beta}$.
We have
	\begin{eqnarray}
		\Sigma_{\e,\beta}( t_{\e}, x_{\e},w_{\e},  y_{\e},v_{\e})
		\geq \Sigma_{\e, \beta}(t^{*}, y^{*}, v^{*}, y^{*}+\beta\sqrt{2\e}\eta_{1},v^{*}+\beta\sqrt{2\e}\eta_{2}),
	\end{eqnarray}
which implies that
	\begin{eqnarray}
	\label{91701}		
	&&\Big(\frac{x_{\e}-y_{ \e}}{\sqrt{2\e}}+\beta\eta_{2}\Big)^{2}+\Big(\frac{w_{\e}-v_{\e}}{\sqrt{2\e}}
		+\beta\eta_{3}\Big)^{2}+ \beta [(t_{\e}-t^{*})^{2}+(x_{\e}-y^{*})^{2}+(w_{\e}-v^{*})^{2}]\nonumber\\
		&\le&\underline{u}(t_{\e},x_{\e},w_{\e})-\bar{u}^{\rho,\theta,\vsi}(t_{\e},y_{\e},v_{\e})
		-\underline{u}(t^{*},y^{*},v^{*})
		+\bar{u}^{\rho,\theta,\vsi}(t^{*}, y^{*}+\beta\sqrt{2\e}\eta_{1},v^{*}+\beta\sqrt{2\e}\eta_{2})\nonumber\\
	&\le&	\frac{2{M}(1+\rho)}{c}+\frac{\theta(T-t^{*}+\varsigma)}{t^{*}}.
		\end{eqnarray}
It follows that $[(x_{\e}-y_{\e})^{2}+(w_{\e}-v_{\e})^{2}]/\e\le C_\beta$ for some constant $C_\beta>0$. Thus, possibly along a subsequence, we have
	$\lim_{\e\rightarrow 0}[(x_{\e}-y_{\e})^{2}+(v_{\e}-w_{\e})^{2}]= 0$.
	By the continuity of the functions $\underline{u}$ and $\bar{u}^{\rho, \th, \vsi}$ and the definition of $(t^*, y^*, v^*)$ we
have
$$\lim_{\e\rightarrow 0}[\underline{u}(t_{\e},x_{\e},w_{\e})-\bar{u}^{\rho,\theta,\vsi}(t_{\e},y_{\e},v_{\e})]\leq M_b=\lim_{\e\to0}
[\underline{u}(t^{*},y^{*},v^{*})
		-\bar{u}^{\rho,\theta,\vsi}(t^{*}, y^{*}+\beta\sqrt{2\e}\eta_{1},v^{*}+\beta\sqrt{2\e}\eta_{2})].
		$$
Therefore, sending $\e \rightarrow 0$ in (\ref{91701})
we obtain that
	\begin{eqnarray*}
		\lim_{\e \rightarrow 0}\Big[\Big(\frac{x_{\e}-y_{ \e}}{\sqrt{2\e}}+\beta\eta_{1}\Big)^{2}+\Big(\frac{w_{\e}-v_{\e}}{\sqrt{2\e}}
		+\beta\eta_{2}\Big)^{2}+ \beta [(t_{\e}-t^{*})^{2}+(x_{\e}-y^{*})^{2}+(w_{\e}-v^{*})^{2}]\Big]\leq 0.
	\end{eqnarray*}
Consequently, we conclude that 	
\begin{eqnarray}
\label{etostar}
		\left\{\ba{lll}
	\dis \lim_{\e\to0}(t_{\e},x_{\e},w_{\e})=\lim_{\e\to0}(t_{\e},y_{\e},v_{\e}) =(t^{*},y^{*},v^{*}), \\
	\dis 	\lim_{\e \rightarrow 0}\Big(\frac{1}{\sqrt{2\e}}(x_{\e}-y_{ \e})+\beta\eta_{1}\Big)^{2}+\Big(\frac{1}{\sqrt{2\e}}
	(w_{\e}-v_{\e})	+\beta\eta_{2}\Big)^{2}= 0.
		\ea\right.
	\end{eqnarray}
In other words, we have shown that
	$$y_{\e}=x_{\e}+\beta \sqrt{2\e}\eta_{1}+o(\sqrt{2\e}), \qq v_{\e}=w_{\e}+\beta \sqrt{2\e}\eta_{2}+o(\sqrt{2\e}),$$
and it then follows from (\ref{91601}) that  $(t_{\e},y_{\e},v_{\e})\in \sD^0_b$ for $\e>0$ small enough.
Namely, we have now returned to the situation of Case 1, with a slightly different penalty function $\f_{\e,\beta}$.
The rest of the proof will follow a similar line of arguments, which we shall present briefly for completeness.
First we apply \cite[Theorem 8.3]{CIL} again to assert that for
any $\delta>0$, there exist $q,\hat{q} \in\hR$ and $A,B \in\cS^2$ such that
	\begin{eqnarray}
	\label{sub}
		\left\{\ba{lll}
		(q, (\xi^{1}_\e+2\beta(x_{\e}-y^{*}),\xi^{2}_\e+2\beta(w_{\e}-v^{*})), A)\in \bar{{\sP}}_{D}^{1,2,+}\underline{u}(t_{\e},x_{\e},w_{\e})\\
		(\hat{q}, (\xi^{1}, \xi^{2}), B)\in \bar{{\sP}}_{D}^{1,2,-}\bar{u}^{\rho,\theta, \vsi}(t_{\e},y_{\e},v_{\e}),
		\ea\right.
	\end{eqnarray}
where  $q-\hat{q}=2\beta(t_{\e}-t^{*})$,
	$\xi^{1}_\e:=(x_{\e}-y_{\e})/\e+2\beta\eta_{1}/\sqrt{2\e}$, $\xi^{2}_\e:=(w_{\e}-v_{\e})/\e+2\beta\eta_{2}/\sqrt{2\e}$, and
	\begin{eqnarray}
		\begin{pmatrix}
			(2 \beta+\frac{1}{\varepsilon})I & -\frac{1}{\varepsilon}I\\
			-\frac{1}{\varepsilon}I & \frac{1}{\varepsilon}I
		\end{pmatrix}+\delta
		\begin{pmatrix}
			(\frac{2}{\varepsilon^{2}}+4\beta^{2}+\frac{4\beta}{\varepsilon})I & -(\frac{2}{\varepsilon^{2}}+\frac{2 \beta}{\varepsilon})I\\
			-(\frac{2}{\varepsilon^{2}}+\frac{2 \beta}{\varepsilon})I & \frac{2}{\varepsilon^{2}}I
		\end{pmatrix}
		\geq
		\begin{pmatrix}
			A & 0\\
			0 & -B
		\end{pmatrix}.			
	\end{eqnarray}
 Now, setting $\delta=\varepsilon$ we have
	\begin{eqnarray}
		\frac{3}{\varepsilon}\begin{pmatrix}
			I & -I\\
			-I & I
		\end{pmatrix}+
		\begin{pmatrix}
			(6\beta+4 \beta^{2}\varepsilon)I & -2 \beta I\\
			-2\beta I & 0
		\end{pmatrix}
		\geq
		\begin{pmatrix}
			A & 0\\
			0 & -B
		\end{pmatrix},
		\label{90901}
	\end{eqnarray}
which implies, in particular,
\bea
\label{90904}
		A_{11}x_{\e}^{2}-B_{11}y_{\e}^{2}\leq \frac{3}{\e}(x_{\e}-y_{\e})^{2}+(6\beta+4 \beta^{2}\e)x^{2}_{\e}-4\beta x_{\e}y_{\e}.
	\eea
Again, as in Case 1 we can easily argue that, without loss of generality,
one may assume that
%
$(q, (\xi^1_\e+2\beta(x_{\e}-y^*), \xi^2_\e+2\beta(w_{\e}-v^*)), A_{11})\in\sP^{+(1,2,1)}_{\sD^*}\underline{u}(t_\e,x_\e,w_\e)$ and
 $(\hat q, (\xi^1_\e,\xi^2_\e), B_{11})\in \sP^{-(1,2,1)}_{\sD} \bar{u}^{\rho,\theta,\vsi}(t_\e,x_\e,w_\e)$).
It is important to notice that, while $(t_\e, y_\e, v_\e)\in \sD^0_b$, it is possible that the point $(t_\e, x_\e, w_\e)$ is on the
boundary of $\sD^*$. Thus it is crucial that viscosity (subsolution) property is satisfied on $\sD^*$, including the boundary
points. Thus,
by Definition \ref{VD2} we have
\beaa
&&	q+\sup_{\gamma\in [0,1], a\in [0,M]}H(t_{\e},x_{\e},w_{\e},\underline{u},\xi^1_\e+2\beta(x_{\e}-y^{*}),\xi_\e^{2}+2\beta(w_{\e}-v^{*}), A_{11}, I[\underline{u}],\gamma,a)\geq 0, \\
&&	\hat{q}+\sup_{\gamma\in [0,1], a\in [0,M]}H(t_{\e},y_{\e},v_{\e}, \bar{u}^{\rho,\theta,\vsi},\xi_\e^{1}, \xi_\e^{2}, B_{11}, I[\bar{u}^{\rho, \th, \vsi}], \gamma,a)\leq0.
\eeaa
Now if we take $(\gamma_{\e},a_{\e})\in {\rm argmax} H(t_{\e},x_{\e},w_{\e}, \underline{u},\xi_\e^{1}+2\beta(x_{\e}-y^{*}),\xi_\e^{2}+2\beta(w_{\e}-v^{*}), A_{11}, I[\underline{u}], \gamma,a)$, then we have
	\begin{eqnarray*}
	0&\le & (q-\hat{q})+	H(t_{\e},x_{\e},w_{\e}, \underline{u},(\xi_\e^{1}+2\beta(x_{\e}-y^{*}),\xi_\e^{2}+2\beta(w_{\e}-v^{*})), A_{11}, I[\underline{u}],\gamma_{\e},a_{\e})\\
	&&-H(t_{\e},y_{\e},v_{\e}, \bar{u}^{\rho,\theta,\vsi}, (\xi_\e^{1}, \xi_\e^{2}), B_{11}, I[\bar{u}^{\rho, \th, \vsi}], \gamma_{\e},a_{\e}),
	\end{eqnarray*}
or equivalently, 
	\begin{eqnarray}
	\label{100302}
		&&(c+\l(w_\e))\underline{u}(t_{\e}, x_{\e},w_{\e})-(c+\l(v_\e))\bar{u}^{\rho,\theta,\vsi}(t_{\e}, y_{\e},v_{\e})\nonumber\\
		&\leq& \frac{1}{2}\sigma^{2}{\gamma_{\e}}^{2}(A^{11}x_{\varepsilon}^{2}-B^{11}y_{\e}^{2})
		+r(x_{\e}-y_{\e})^{2}/\e+2(x_{\e}-y_{\e})r\beta \eta_{1}/\sqrt{2\e} \nonumber \\
		&&+2\beta[(r x_{\e}+p-a)(x_{\e}-y^{*})+(w_{\e}-v^{*})]
		+2\beta(t_{\e}-t^{*})\nonumber\\
		&&+\l(w_{\e})\int_{0}^{x_{\e}} \underline{u}(t_{\e},x_{\e}-u,0)dG(u)-\l(v_{\e})\int_{0}^{y_{\e}} \bar{u}^{\rho,\theta,\vsi}(t_{\e},y_{\e}-u,0)dG(u)\\
		&\leq& (3\sigma^{2}{\gamma_{\e}}^{2}/2
		+r)(x_{\e}-y_{\e})^{2}/\e+2(x_{\e}-y_{\e})r\beta \eta_{1}/\sqrt{2\e}  \nonumber \\
		&&+2\beta[(r x_{\e}+p-a)(x_{\e}-y^{*})+(3+2 \beta\e)x^{2}_{\e}-2 x_{\e}y_{\e}+(w_{\e}-v^{*})+(t_{\e}-t^{*})]\nonumber \\
		&&+\l(w_{\e})\int_{0}^{x_{\varepsilon}} \underline{u}(t_{\e},x_{\e}-u,0)dG(u)-\l(v_{\e})\int_{0}^{y_{\e}} \bar{u}^{\rho,\theta,{\color{purple2}\vsi}}(t_{\e},y_{\e}-u,0)dG(u).\nonumber
			\end{eqnarray}
	First sending $\varepsilon\rightarrow 0$ then sending $\beta\rightarrow 0$, and noting (\ref{etostar}),
	we obtain from (\ref{100302}) that
	\begin{eqnarray*}
		(c+\l(v^{*}))M_b\leq
		\l(v^{*})(\int_{0}^{y^{*}}( \underline{u}(t^{*},y^{*}-u,0)- \bar{u}^{\rho,\theta,\vsi}(t^{*},y^{*}-u,0))dG(u))
		\leq \l(v^{*}) M_b.
	\end{eqnarray*}
Again, this is a contradiction as $c>0$ and $M_b>0$. The proof is now complete.
\qed

%
%
%
%
%

\end{document}